\documentclass [11pt,oneside,a4paper,mathscr]{amsart}

\usepackage{amscd,amsmath,amssymb,euscript}
\usepackage[frame,cmtip,curve,arrow,matrix,line,graph]{xy}
\usepackage{graphicx}

\title{Helices on del Pezzo surfaces and tilting Calabi-Yau algebras}
\author{Tom Bridgeland and David Stern}

\date{}

\newtheorem{thm}{Theorem}[section]

\newtheorem{cor}[thm]{Corollary}
\newtheorem{prop}[thm]{Proposition}
\newtheorem{lemma}[thm]{Lemma}
\newenvironment{pf}{\paragraph{Proof}}{\qed\par\medskip}
\theoremstyle{definition}
\newtheorem{defn}[thm]{Definition}
\newtheorem{remark}[thm]{Remark}
\newtheorem{remarks}[thm]{Remarks}
\newtheorem{example}[thm]{Example}
\newtheorem{examples}[thm]{Examples}

\renewcommand{\leq}{\leqslant}
\renewcommand{\geq}{\geqslant}

\renewcommand{\H}{\mathbb{H}}
\newcommand{\E}{\mathbb{E}}
\newcommand{\Mod}{\operatorname{Mod}}
\newcommand{\End}{\operatorname{End}}
\newcommand{\K}{{\operatorname{K}}}

\newcommand{\D}{{\operatorname{D}}}
\newcommand{\Coh}{\operatorname{Coh}}

\newcommand{\isom}{\cong}
\renewcommand{\L}{{L}}
\newcommand{\tensor}{\otimes}

\newcommand{\PP}{\operatorname{\mathbb P}}
\newcommand{\C}{\mathbb C}

\newcommand{\blob}{{\scriptscriptstyle\bullet}}
\newcommand{\F}{\mathbb F}
\newcommand{\G}{\mathbb G}

\newcommand{\Z}{\mathbb Z}
\newcommand{\A}{\mathcal A}
\newcommand{\B}{\mathcal B}

\newcommand{\CC}{\mathcal C}
\newcommand{\OO}{\mathcal O}
\newcommand{\into}{\hookrightarrow}
\newcommand{\id}{\operatorname{id}}
\newcommand{\Ext}{\operatorname{Ext}}
\newcommand{\Hom}{\operatorname{Hom}}

\newcommand{\lRa}[1]{\xrightarrow{\ #1\ }}
\newcommand{\lra}{\longrightarrow}
\newcommand{\R}{{R}}

\newcommand{\ec}{{exceptional collection }}
\newcommand{\ecs}{exceptional collections }
\renewcommand{\O}{\OO}
\newcommand{\DQcoh}{\D\operatorname{Qcoh}}
\newcommand{\op}{{\operatorname{op}}}

\newcommand{\CQ}{\C Q}
\newcommand{\CY}{$\text{CY}_3$ }

\renewcommand{\S}{\mathbb{S}}
\newcommand{\ideal}{\vartriangleleft}
\newcommand{\Proj}{\operatorname{Proj}}
\newcommand{\fin}{\operatorname{fin}}

\begin{document}

\begin{abstract}We study tilting for a class of Calabi-Yau
algebras associated to  helices on Fano varieties. We do this by relating
the tilting operation to
mutations of exceptional collections. For helices on
del Pezzo surfaces the  algebras are of dimension
three, and using an argument of Herzog, together with results of
Kuleshov and Orlov, we obtain a complete description of the
tilting process in terms of quiver mutations.\end{abstract}

\maketitle

\section{Introduction}
 In the context of homological algebra, tilting
is a fundamental construction that relates neighbouring
t-structures in a triangulated category. It first appears at this
level of generality in a paper of Happel, Reiten and Smal{\o}
\cite{hrs}. The basic idea however goes back much further, and the
name was introduced by Brenner and Butler \cite{bb}, who studied
the process in the context of categories of
representations of quivers. 
More recently  tilting  for three-dimensional
Calabi-Yau algebras has been  related to cluster mutations \cite{k,kd,ks}, and to Seiberg duality in theoretical physics \cite{bd}.
\smallskip

The aim of this paper is to study tilting for a class of
 Calabi-Yau algebras associated to helices of coherent sheaves on Fano varieties. Our basic
tool will be the theory of exceptional collections, developed in
the Rudakov seminar \cite{r}. Before stating our main result
(Theorem \ref{main}) we will recall the definition of tilting in
the context of algebras defined by quivers with relations, and
introduce the class of algebras to be studied.

\subsection*{Notation}
All algebras and varieties will be over the complex numbers. All
modules will be right modules. We write $\D(X)$ for the bounded
derived category of coherent sheaves on a variety $X$, and $\D(A)$
for the homotopy category of bounded complexes of
finitely-generated projective modules over an algebra $A$. If $A$
is noetherian and of finite global dimension  this is equivalent
to the bounded derived category of finitely-generated $A$-modules.
Finally, we write $\D_{\fin}(A)\subset\D(A)$ for the subcategory
of complexes with finite-dimensional cohomolgy.

\subsection{Quivers and tilting}
\label{saville}

 Let $Q$ be a quiver specified in the usual way
by a set of vertices $Q_0$, a set of arrows $Q_1$, and source and
target maps $s,t\colon Q_1\to Q_0$. We shall always assume that the number of arrows  $n_{ij}$ from vertex $i$ to vertex $j$ is finite.
The path algebra $\CQ=\bigoplus _{l\geq 0} \CQ_l$ is
graded by path length, and the degree zero part $\S=\CQ_0$ is a semisimple
 ring with a basis of orthogonal idempotents $e_i$ indexed by the vertices $Q_0$.

  Given a  two-sided ideal $I\ideal\CQ$ generated by linear
combinations of paths of length at least two, let
\[A=A(Q,I)=\CQ/I\] denote the corresponding quotient algebra.
Then $A$ is an augmented  $\S$-algebra
with augmentation ideal
 $A_+\ideal A$ spanned by paths of positive length.
 The underlying quiver
  $Q$ is determined by the augmented algebra $A$, since
 \begin{equation}
 \label{nearlydone}
 n_{ij} =\dim_{\C} e_j \cdot\big( A_+/(A_+)^2\big)\cdot
 e_i.\end{equation}
 We shall
refer to augmented  algebras of this form as \emph{quiver
algebras}.

Suppose that $A=A(Q,I)$ is a quiver algebra.
To each vertex $i\in Q_0$ there corresponds an indecomposable projective
module $P_i=e_i A$, and a one-dimensional simple module $S_i$, on which all elements of $A_+$ act by zero. It is easy to see that
\begin{equation}\label{mortlock}n_{ji}=\dim_\C \Ext^1_A(S_i,S_j).\end{equation}
There is also a canonical map
\begin{equation*} \bigoplus_{j\in Q_0} P_j^{\oplus n_{ji}} \lra P_i\end{equation*}
which, viewed as a complex of modules concentrated in degree $0$ and $1$,
defines an object $R_i$ of the category $K^b\Proj(A)$.

\begin{defn}\label{tommy}
We say that  quiver algebras $A=A(Q,I)$ and $A'=A(Q',I')$  are related by a
\emph{vertex tilt} if  there is an equivalence of categories
\[\Psi \colon \D(A) \to \D(A'),\] a bijection $\psi\colon Q_0\to Q'_0$,
and a vertex $i\in Q_0$ such that
\[\Psi(P_j)=P'_{\psi(j)}\text{ for }j\neq i\quad \text{and}\quad\Psi(R_i)=P'_{\psi(i)}.\]
More precisely we say that $A'$ is the left tilt of $A$ at the vertex $i$, and  that $A$ is
the right tilt of $A'$ at the vertex $\psi(i)$.
\end{defn}



\begin{remark}
\label{bandy} Suppose that the quivers underlying the algebras $A$
and $A'$ have no loops. The equivalence $\Psi$ restricts to an
equivalence
\[\Psi\colon \D_{\fin}(A)\to \D_{\fin}(A').\]Define objects
$U_j\in \D_{\fin}(A)$ for $j\in Q_0$ by the relation
$\Psi(U_j)=S_{\psi(j)}$.  Then   $U_i=S_i[-1]$, whereas for $j\neq
i$ the object $U_j$ is the universal extension
\begin{equation} \label{first} 0\lra S_j\lra U_j\lra
\Ext^1_{A}(S_i,S_j)\tensor {S_i} \lra 0.\end{equation} For a proof
of these facts see Lemma \ref{golden}. It is an easy consequence
of this that when the algebras $A$ and $A'$ are noetherian of
finite global dimension the inverse image under $\Psi$ of the
standard t-structure on $\D_{\fin}(A')$ is related to the standard
t-structure on $\D_{\fin}(A)$  by an abstract tilt in the sense of
\cite{hrs}. The relevant torsion pair  has torsion part consisting
of direct sums of the module $S_i$.
\end{remark}

\subsection{Tilting \CY quiver algebras}
We shall say that a quiver algebra $A=A(Q,I)$ is (weakly) Calabi-Yau of dimension $d$, or just CY$_d$, if $A$ has global dimension $d$
and the shift functor $[d]$ is a Serre functor on $\D_{\fin}(A)$. This last condition means that there are functorial isomorphisms
\[\Hom_A(M,N)\isom \Hom_A(N,M[d])^*\]
for all objects $M,N\in\D_{\fin}(A)$.
We shall be particularly interested in the case of \CY quiver algebras.
The  combinatorics of the tilting process for such algebras
 can be described by a  rule known as quiver mutation, as
we now explain.

The \CY condition implies that the Euler matrix
\[\chi(i,j)=\chi(S_i,S_j)=\sum_{i=0}^3 (-1)^i \dim_{\C} \Ext^i_A (S_i,S_j).\]
 is  skew-symmetric.
Suppose that the underlying quiver $Q$ has no loops or
oriented 2-cycles. Then by \eqref{mortlock}
\[n_{ij} =  \Biggl\{\begin{array}{ll} \chi({i,j}) & \text{ if }\chi({i,j})>0, \\ 0 &\text{
otherwise.}\end{array}\Biggr.\]
Thus the structure of the quiver is completely determined by the
Euler matrix.

\begin{lemma}
\label{sun} Suppose $A=A(Q,I)$ and $A'=A(Q',I')$ are \CY quiver
algebras related by a tilt at the vertex $i\in Q_0$ as in
Definition \ref{tommy}.  Then the Euler form for the algebra $A'$
is
\[\chi(\psi(j),\psi(k))= \biggl\{\begin{array}{ll} -\chi(j,k) & \text{ if }i\in \{j,k\},\\ \chi(j,k)&\text{
if } i\notin \{j,k\}\text{ and }\chi(i,j)\, \chi({i,k})\geq 0, \\
\chi(j,k)+|\chi(i,j)|\cdot \chi(i,k)  &\text{ if } i\notin \{j,k\}\text{
and }\chi(i,j)\, \chi(i,k)\leq 0,\end{array}\biggr.\]
\end{lemma}

\begin{pf}
This follows easily from \eqref{first} using the additivity
property of $\chi$ on exact sequences. \end{pf}

Thus if we also assume that $Q'$ has no loops or 2-cycles  the
quiver $Q'$  is completely determined by the quiver $Q$ and the
vertex $i$. The resulting transformation law $\mu_i\colon Q\mapsto
Q'$ is called \emph{quiver mutation} and is easily checked to be
an involution. It is very important to note however that one does
not usually know \emph{a priori} that the quiver $Q'$ has no
2-cycles, so that without further information, one cannot
guarantee that the quiver underlying the algebra $A'$ is indeed
given by the above rule.

\begin{remark}
\label{twits} Keller and Yang \cite{kd} have recently gone further
and described how the ideals of relations $I$ and $I'$ are related
under a tilt. Take assumptions as in Lemma \ref{sun}. For
simplicity let us also assume that $A$ and $A'$ are graded. Work of
Bocklandt \cite[Theorem 3.1]{bok} then shows that the relations
 can be encoded in compact form in a
potential. Thus we can write
\[A=A(Q,W)=\CQ/(\partial_a W:a\in Q_1)\]
for some non-uniquely defined element $W\in \CQ/[\CQ,\CQ]$.
We can  similarly write $A'=A(Q',W')$. Then Keller and Yang \cite[Theorem 3.2]{kd} (see also \cite[Theorem 9.2]{k})
show that the potentials $W$ and $W'$ can be chosen  so that $(Q',W')$ is obtained from $(Q,W)$ by
a simple combinatorial rule described explicitly by
 Derksen, Weyman and Zelevinsky \cite{DWZ}.
\end{remark}


\subsection{Rolled-up helix algebras}

Let $Z$ be a smooth projective  variety and let $\omega_Z$ denote
its canonical line bundle. In the cases of most interest to us $Z$ will be a Fano variety, which is to say that the dual of $\omega_Z$ is ample.

\begin{defn}
A \emph{helix of sheaves} on $Z$ of period $n$ is an infinite  collection of coherent sheaves 
$\H=(E_i)_{i\in \Z}$ such that for all $i\in \Z$ one has

\begin{itemize} \item[(i)] $(E_{i+1}, \cdots, E_{i+n})$ is a full exceptional
collection, \smallskip \item[(ii)] $E_{i-n}=E_i \tensor \omega_Z.$

\end{itemize}
\end{defn}
We recall many aspects of the theory of exceptional collections and helices in
Sections \ref{exxc} and \ref{heli} below.
Associated to a helix $\H$ is a graded algebra
\[A(\H)={\bigoplus_{k\geq 0}}\prod_{j-i=k} \Hom_{Z} (E_i,E_j)\]
equipped with a $\Z$-action induced by twisting by $\omega_Z$
\[-\tensor \omega_Z\colon \Hom_{Z}(E_i,E_j) \to \Hom_{Z}(E_{i-n},E_{j-n}).\]
The rolled-up helix algebra $B(\H)$ is defined to be the
subalgebra of $A(\H)$ consisting of elements invariant under this
action.

 A helix $\H=(E_i)_{i\in \Z}$ is called
\emph{geometric} if for all $i<j$ one has the additional condition
\[\Hom_{Z}^k(E_i,E_j)=0\text{ unless }k=0.\]
%
The following result was mostly proved in \cite[Prop. 4.1]{br}, see also \cite[Prop. 7.2]{vdb}. For completeness we give
a proof in Section \ref{proof}.


\begin{thm}
\label{snow} Let $B=B(\H)$ be the rolled-up helix algebra of a
geometric helix on a smooth projective Fano variety $Z$ of
dimension $d-1$. Then $B$ is a graded CY$_d$ quiver algebra which
is noetherian and finite over its centre. Moreover, there are
equivalences of categories $\D(B)\isom \D(Y)$
 where $Y$ is the  total space of the line bundle $\omega_Z$.
\end{thm}

%
The tilting operation for rolled-up helix algebras was studied in
\cite{br} under the additional assumption
\begin{equation}\label{dj}\operatorname{rank} K(Z)=\dim Z +1.\end{equation}
This case is much simpler than the general one. In particular, the
quiver underlying $B(\H)$  is always an oriented $n$-gon with
various numbers of arrows along the edges, and the combinatorics
of the tilting process  can be described explicitly using an
action of the affine braid group.
\smallskip

In this paper we shall be most interested in  the case when $Z$ is
a del Pezzo surface, that is a Fano variety of dimension
2. Recall that any such surface is isomorphic either to\begin{itemize}
\item[(a)]
$\PP^1\times\PP^1$, or
\item[(b)] the blow-up of $\PP^2$ at $m \leq 8$
points.
\end{itemize}
The only case satisfying \eqref{dj} is $Z=\PP^2$.

Thanks to work of Karpov and Nogin \cite{kn} it is known that
geometric helices exist on all del Pezzo surfaces (see Example
\ref{thisistheend} below). The main content of the following
result is that it provides an interesting class of \CY algebras,
closed under vertex tilts, where the quiver mutation rules really
apply, in that one knows \emph{a priori} that there are no loops
or 2-cycles.

\begin{thm}
\label{main} Let $B=B(\H)$ be the rolled-up helix algebra of a
geometric helix $\H$ on a del Pezzo surface $Z$. Then $B$ is a
graded \CY quiver algebra whose underlying quiver $Q$ has no loops
or oriented 2-cycles. For any vertex $i$ of $Q$ there is another
geometric helix $\H'$ on $Z$ such that the algebra $B(\H')$ is the
(left and right) tilt of $B(\H)$ at the vertex $i$.
\end{thm}

It follows from Theorem \ref{main} and Remark \ref{twits}  that
tilting for rolled-up helix algebras on del Pezzo surfaces is
completely described by the combinatorial quiver mutation process of \cite{DWZ}.
It is worth noting that the same combinatorics was considered earlier by physicists studying Seiberg duality for quiver gauge theories, see for example  \cite{Hanany}.

%
%
%
%
%



%
\subsection{An example}
As an example take $Z=\PP^1\times \PP^1$ with its projections
$\pi_1,\pi_2\colon Z\to \PP^1$. We use the standard notation
\[\O(a,b)=\pi_1^* \O_{\PP^1}(a) \tensor \pi_2^* \O_{\PP^1}(b).\]
 The canonical bundle of $Z$ is $\O(-2,-2)$.
The sequence
\[\cdots ,\O,\O(1,0),\O(0,1),\O(1,1),\O(2,2),\cdots\]
is a geometric helix of  period 4. The corresponding rolled-up
helix algebra has quiver
\[
\xymatrix@C=5em{ \bullet\ar[d]_{2}\ar[r]^{2} & \bullet\ar[d]^{2}\\
\bullet\ar[r]_{2} &\bullet \ar[ul]^{4} }
\]
Tilting at the top right vertex gives another rolled-up helix
algebra with quiver
\[
\xymatrix@C=5em{ \bullet\ar[r]^{2} & \bullet\ar[d]^{2}\\
\bullet\ar[u]^{2} &\bullet \ar[l]^{2} }
\]
corresponding to the geometric helix
\[\cdots ,\O,\O(1,0),\O(1,1),\O(1,2),\O(2,2),\cdots\]
Continuing the tilting process one obtains an infinite web of \CY
algebras, part of which is shown in Figure 1, with the associated
quivers. \smallskip



\begin{figure}[h]
\includegraphics[width=.5\linewidth]{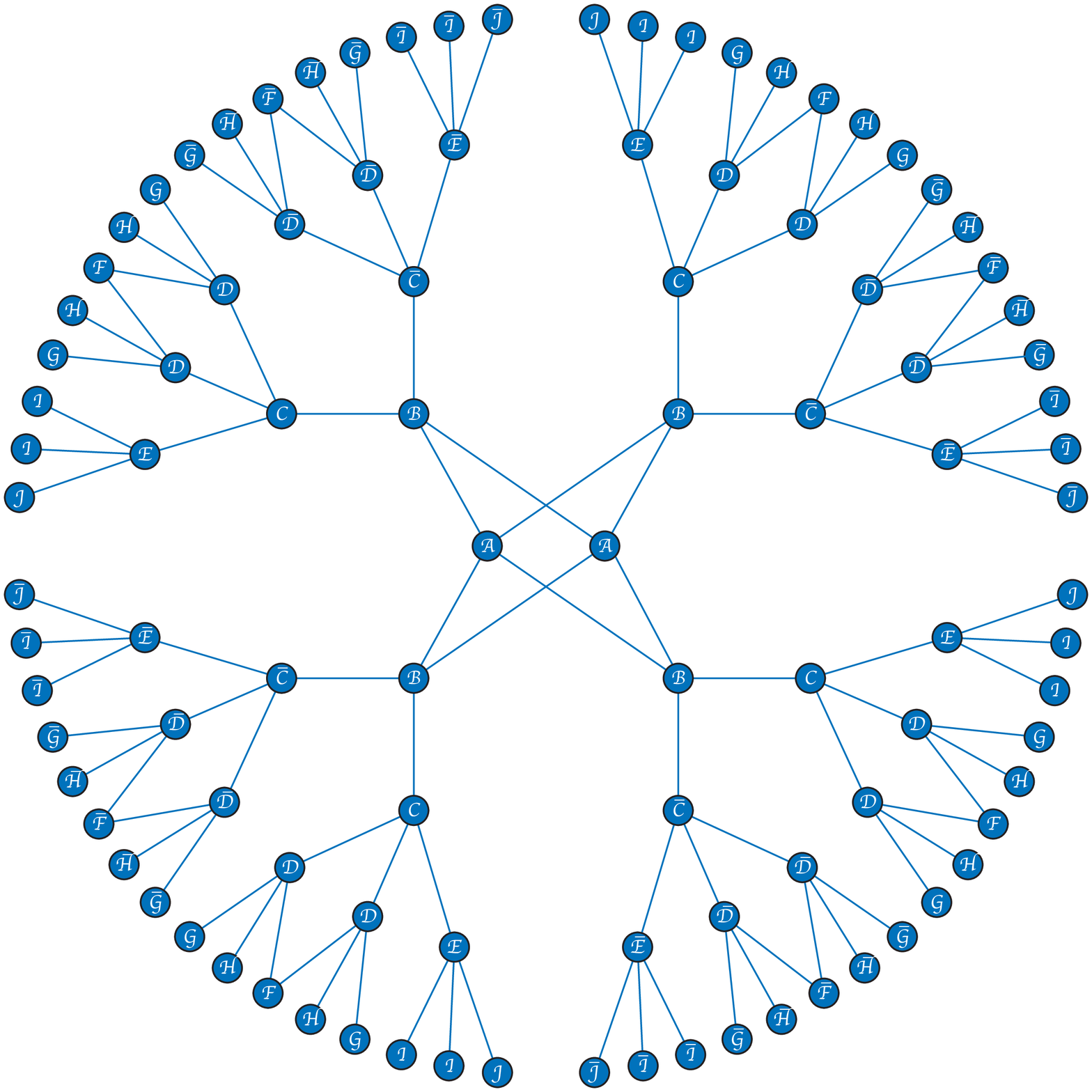}\\
\includegraphics[width=0.65\linewidth]{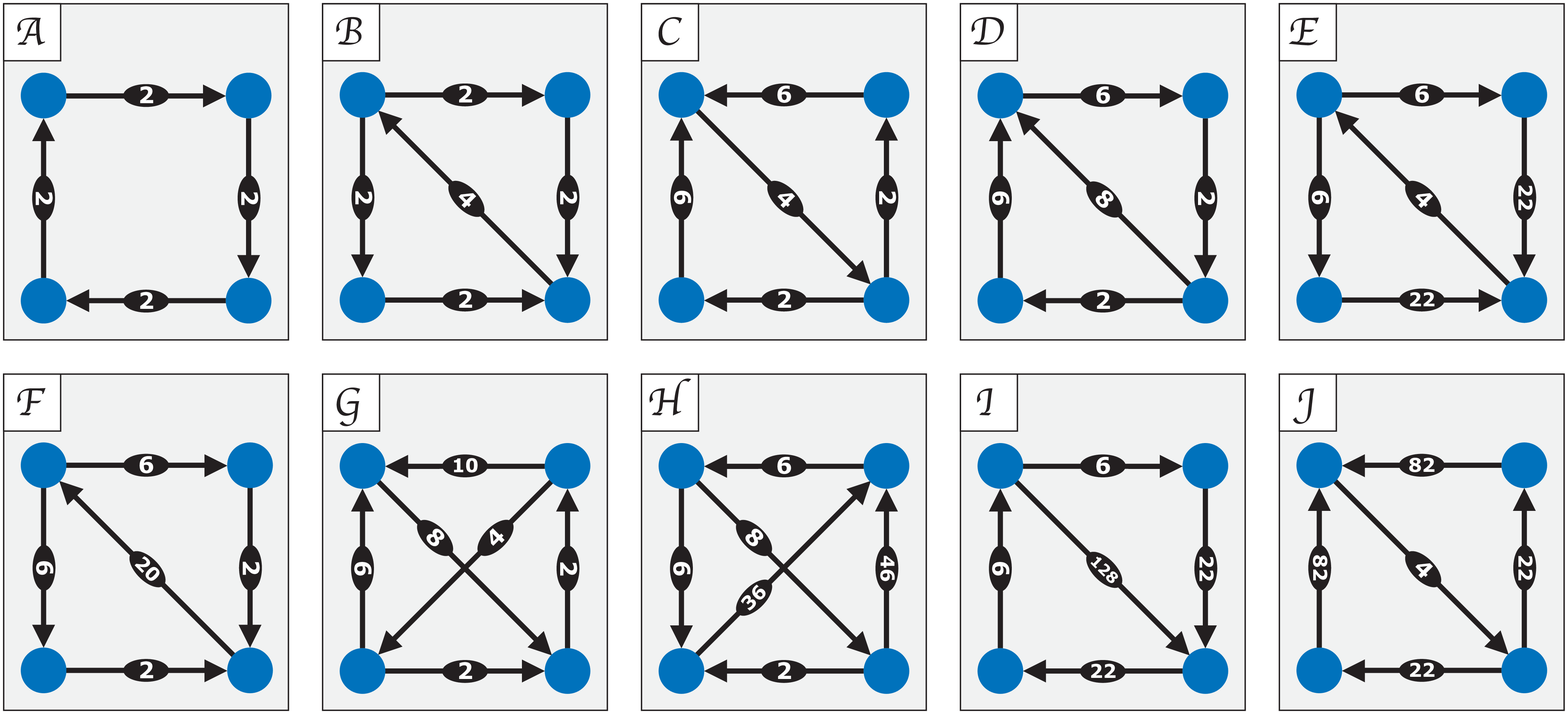}
\caption{ A web of Calabi-Yau algebras.}
\end{figure}


\subsection*{Acknowledgements}
This paper is a partial write-up of the second authors PhD thesis.
It has benefited from useful conversations with Bernhard Keller,
Alastair King and Michel van den Bergh. The authors are also very
grateful to Ben Eldridge for help with producing Figure 1.


\section{Exceptional collections and mutation functors}
\label{exxc}

In this section we review the basic definitions and results
concerning exceptional collections and mutations following
\cite{bo,bk}. This material is by-now standard,  but we have
departed slightly from the usual convention in defining mutation
functors using the natural categorical shifts. This eliminates
unnecessary shifts from several formulae.

\subsection{Assumptions} Throughout $\D$ will be a fixed $\C$-linear
triangulated category. We always assume that
\begin{itemize}
\item{} $\D$ is of \emph{finite type}, i.e. for any two objects
$A,B$ the vector space \[\bigoplus_{i\in \Z} \Hom^i_\D(A,B) \] is
finite-dimensional.
\smallskip
\item{} $\D$ is \emph{algebraic} in the sense of Keller \cite[Section 3.6]{kedg}.
\end{itemize}

The main examples we have in mind are
\begin{itemize} \item[(a)] The bounded derived category of
coherent sheaves $\D(Z)$ on a smooth projective variety $Z$. \smallskip
 \item[(b)] The bounded
derived category of finitely generated $A$-modules $D(A)$ over a
finite-dimensional algebra $A$ of finite global dimension.
\end{itemize}

It will usually also be the case that

\begin{itemize}
 \item{}  $\D$ is \emph{saturated} \cite{bk}, i.e.
 all triangulated functors
\[\D\to \D(\C), \quad \D^\op\to \D(\C),\]
are representable.
\end{itemize}

Note that if $\D$ has a full exceptional collection then this
condition is automatic (see \cite[Corollary to Theorem 2.10]{bk}).
Moreover the two classes of examples (a) and (b) above are saturated  by
\cite[Theorem 1.1]{bvdb} and \cite[Theorem 2.11]{bk} respectively.

\subsection{Exceptional collections}

An object $E\in \D$ is said to be \emph{exceptional} if
\[\Hom_{\D}^k(E,E)=\biggl\{\begin{array}{ll} \C & \text{ if }k=0, \\ 0 &\text{
otherwise.}\end{array}\biggr.\]
An \emph{exceptional collection}
$\E\subset\D$
is a sequence of exceptional objects \[\E=(E_1,\cdots
,E_{n})\]  such that if $1\leq i<j\leq n$ then
$\Hom^\blob_{\D}(E_j,E_i)=0$.
\smallskip

Given an \ec $\E$ the right orthogonal subcategory to $\E$ is
the full triangulated subcategory
\[
\E^{\perp} =\{X\in \D : \Hom_{\D}^{\blob} (E,X)=0 \text{ for }
E\in\E\}.\]
Similarly, the left orthogonal subcategory to $\E$ is
\[{^\perp} \E = \{X\in \D: \Hom_{\D}^{\blob}(X,E)=0 \text{ for } E\in \E\}.\]
When $\E$ consists of a single object $E$ we just write $E^\perp$ and ${^\perp} E$.
\smallskip

We write $\langle \E\rangle\subset \D$ for the smallest full
triangulated subcategory of $\D$ containing the elements of an \ec
$\E\subset \D$. An exceptional collection $\E\subset\D$  is said to be \emph{full}
if $\langle \E\rangle = \D$. By Lemma \ref{full} this is
equivalent to assuming either that $\E^{\perp}=0$ or that
${^\perp}\E=0$.

\begin{example}
\label{bib} We give examples for  $\D=\D(Z)$ with $Z$ a projective variety.
\begin{itemize}
\item[(a)] The sequence
\[\big(\O,\O(1),\cdots, \O(k)\big)\]
is a full exceptional collection on $Z=\PP^{k}$.

\item[(b)] The sequence
\[\big(\O,\O(1,0),\O(0,1),\O(1,1)\big)\]
is a full exceptional collection on $Z=\PP^1\times \PP^1$.
\item[(c)] The sequence
\[\big(\O,\O(1,0),\O(-2,1),\O(-1,1)\big)\]
is another full exceptional collection on $Z=\PP^1\times \PP^1$.
\end{itemize}
\end{example}

It is easy to see (see Lemma \ref{grot}) that the classes of the
elements of a full exceptional collection  form a basis for the
Grothendieck group $\operatorname{K}(\D)$. Thus the length of a
full \ec (if one exists) is an invariant of the category $\D$.

\subsection{Mutation functors}
 \label{mut1}

 Suppose $E\in \D$ is
a exceptional. Given an object $X\in {^\perp E}$  the
\emph{left mutation of  $X$ through $E$} is the object $L_E (X)\in
E^\perp$ defined up to isomorphism by the triangle
\begin{equation} \label{left} \Hom_{\D}^{\blob}(E,X)\tensor
E\lRa{ev} X\lra \L_E (X) ,\end{equation} where $ev$ denotes the
evaluation map. Similarly, given $X\in E^\perp$, the \emph{right
mutation of $X$ through $E$} is the object $\R_E (X)\in {^\perp E}$
defined by the triangle
\begin{equation} \label{right}
\R_E (X)\lra X \lRa{coev}\Hom_{\D}^{\blob}(X,E)^*\tensor
E,\end{equation} where $coev$ denotes the coevaluation map. It is
easy to check that these two operations define mutually inverse
equivalences of categories
\[\xymatrix{ {^\perp E} \ar@/^1pc/[rr]^{\L_E} && E^\perp \ar@/^1pc/[ll]^{\R_E}}\]

 We also consider mutations through
exceptional collections.
 Suppose $\E=(E_1,\cdots, E_{k})$ is an exceptional collection in $\D$. Define the
left mutation of an object $X\in {^\perp \E}$ through the
collection $\E$ to be the object
\[\L_{\E}(X) = \L_{E_{1}}\cdots \L_{E_{k}}(X)\in \E^\perp.\]
Similarly, define the right mutation of $X\in \E^\perp$ through
the collection $\E$
 to be the object
\[\R_{\E}(X) = \R_{E_{k}}\cdots \R_{E_{1}}(X)\in {^\perp \E}.\]
Once again these two operations define mutually inverse functors
\[\xymatrix{ {^\perp \E} \ar@/^1pc/[rr]^{\L_\E} && \E^\perp \ar@/^1pc/[ll]^{\R_\E}}\]

 Note that passing from the category $\D$ to its
opposite category $\D^\op$ exchanges  left and right orthogonal
subcategories and also left and right mutation functors. This
symmetry between left and right  will be a constant feature in
this paper, and we will often make statements just for left
mutations, safe in the knowledge that the corresponding statements
for right mutations can be deduced by passing to $\D^\op$.
\begin{remarks}
\label{gather}
 \begin{itemize}
 \item[(a)] If $(E,F)$
is an exceptional pair, then
\[\langle(L_E(F),E)\rangle=\langle
(E,F)\rangle=\langle(F,R_F(E))\rangle.\] This is easily checked
directly from the definition.
\smallskip

 \item[(b)]
 There is a more
categorical approach to the mutation functors which we briefly
recall in  Appendix \ref{semi}. From this approach one obtains the following characterisation of mutations.
Suppose
\[E\lra X\lra Y\]
is a triangle in $\D$ such that $E\in \langle\E\rangle$, $X\in {^\perp}\E$ and $Y\in\E^\perp$. Then it follows that $Y=\L_\E(X)$ and $X=R_\E(Y)$.
\smallskip

\item[(c)]In particular
it follows from (b) that the functors $\L_\E$ and $R_\E$ depend only on the subcategory
$\langle\E\rangle$ and not on the particular choice of
 exceptional collection $\E$.
\end{itemize}
\end{remarks}

\subsection{Homomorphism algebra}
\label{homo}

Let
$\E=(E_1,\cdots,E_n)$ be an \ec in $\D$.
One can associate to $\E$ a finite-dimensional graded algebra
\[A=A(\E)=\End_{\D}(\bigoplus_{i=1}^n E_i)=\bigoplus_{k\geq 0} \bigoplus_{j-i=k} \Hom_{\D}(E_i,E_j) \] called the
\emph{homomorphism algebra} of  $\E$.
\smallskip

The degree zero part of $A$ is an $n$-dimensional
semisimple algebra.  It follows (see Lemma \ref{simon}) that $A$
can be presented as a quiver algebra $A(Q,I)$ for a unique quiver
$Q$.
The vertices of $Q$ correspond to
the elements of the collection $\E$, and hence are naturally indexed by the set $\{1,\cdots,n\}$. By \eqref{nearlydone} the number of arrows $n_{ij}$ in $Q$ from vertex $i$ to vertex $j$
 is 0 unless $i<j$ in which case $n_{ij}$ is the dimension of the  cokernel of the map
\[\bigoplus_{i<k<j}\Hom_\D(E_i,E_k)\tensor \Hom_\D(E_k,E_j)\lra\Hom_\D(E_i,E_j).\]
The fact that $n_{ij}=0$ unless $i<j$ implies that $A$ is of global dimension $\leq n$ (see
Remark \ref{rich}).

\begin{example}
\label{exmp} Consider the full exceptional collection
\[\E=(\O,\O(1,0),\O(0,1),\O(1,1))\]
from Example \ref{bib}(b). Using the above rule it is easy to see that the quiver underlying $A(\E)$ is
\[
\xymatrix@C=5em{ \bullet\ar[d]_{2}\ar[r]^{2} & \bullet\ar[d]^{2}\\
\bullet\ar[r]_{2} &\bullet 
 }
\]
where the numbers on the edges indicate numbers of arrows. We give another way to compute this quiver in Example \ref{shere} below.
\end{example}

An exceptional collection $\E=(E_1,\cdots,E_n)$ is said to be
\emph{strong} if for all $i,j$
\[\Hom^k_{\D}(E_i,E_j)=0 \text{ unless }k=0.\]
In Example \ref{bib} the collections (a) and (b) are strong, but not
the collection (c).

\begin{thm}\label{RB}
Suppose $\D$ is an algebraic triangulated category of finite type.
Suppose $\E=(E_1,\cdots,E_n)$ is a full, strong \ec in $\D$ and
let $A=A(\E)$ be its homomorphism algebra. Then there is an
equivalence
\[\Phi_{\E}\colon \D(A) \lra \D\]
sending the rank one free module $A$ to the object
$E=\bigoplus_{i=1}^n E_i$.
\end{thm}

\begin{pf}In the case when $\D=\D(Z)$ for a smooth projective variety
$Z$ this result is due to Bondal \cite[Theorem 6.2]{bo} and
Rickard \cite{ri}. As stated it follows from general results of
Keller, see \cite[Theorem
3.8(a)]{kedg}.\end{pf}


\subsection{Dual collection}
The following simple result will be very important later.

\begin{lemma}
\label{dual} Let $\E=(E_1,\cdots, E_n)$ be a full exceptional
collection and define \begin{equation} \label{chee}F_j=\L_{E_1}
\cdots \L_{E_{j-1}} (E_j) \quad 1\leq j \leq n \end{equation} Then
$\F=(F_n,\cdots,F_1)$ is a full exceptional collection and
\begin{equation}
\label{freddie}\Hom^{k}_{\D}(E_i,F_j)=\left\{
\begin{array}{ll}
\C &\text{$i=j$ and $k=0$}, \\
0 &\text{otherwise.}
\end{array}\right.\end{equation}
Conversely, if a sequence of objects $(F_n,\cdots,F_1)$ satisfy
\eqref{freddie} then they are given by \eqref{chee}.
\end{lemma}

\begin{pf}
A special case of this appears in \cite[Lemma 5.6]{bo}. For the convenience of the reader we give
a proof in Appendix B.
\end{pf}

We shall call the collection $\F$ of Lemma \ref{dual} the
\emph{dual collection} to $\E$. The importance of dual collections
is as follows. Suppose $\E$ is a full, strong exceptional
collection and let $A=A(\E)$ denote the corresponding homomorphism
algebra. Let $S_i$ be the one-dimensional simple $A$-module corresponding to vertex $i$, and $P_i=e_i A$ the indecomposable
projective module.

\begin{lemma}
\label{reallydull} Under the equivalence $\Phi_\E$ of Theorem
\ref{RB}
\[\Phi_\E (P_i)=E_i, \qquad \Phi_\E(S_i)=F_i.\]
\end{lemma}

\begin{pf}
The first statement holds because the objects $E_i$ are the
indecomposable summands of the object $E$, and the modules $P_i$
are the indecomposable summands of the module $A$. The second
statement follows from the fact
that\[\Hom^{k}_{A}(P_i,S_j)=\left\{
\begin{array}{ll}
\C &\text{$i=j$ and $k=0$}, \\
0 &\text{otherwise,}
\end{array}\right.\]
together with the uniqueness statement of Lemma \ref{dual}.
\end{pf}

\begin{remark}
\label{extra} Later we shall often use the following fact. In the
situation of Theorem \ref{RB} the standard t-structure on $\D(A)$
corresponds under the equivalence $\Phi_\E$ to a bounded,
non-degenerate t-structure $\D$.
The heart of this t-structure is equivalent to $\Mod(A)$ and is a finite-length abelian category whose simple objects
are the elements of the dual collection $\F$. An element $X\in\D$ lies in this heart precisely if $\Hom_\D^k(E_i,X)=0$ for $1\leq i\leq n$ and $k\neq 0$.
\end{remark}

By \eqref{mortlock} from the introduction and Lemma
\ref{reallydull}, working out the dimensions of the spaces $\Hom^1_\D(F_i,F_j)$ gives another way to compute the quiver
underlying the homomorphism algebra of a full, strong exceptional collection.

\begin{example}
\label{shere} Take $Z=\PP^1\times \PP^1$ and consider again the
full, strong exceptional collection
\[\E=(\O,\O(1,0),\O(0,1),\O(1,1)).\]
The corresponding dual collection $\F$  is
\[\F=\big(\O(-1,-1)[2], \O(0,-1)[1], \O(-1,0)[1],\O\big).\]
Computing the dimensions of the spaces $\Hom^1_Z(F_i,F_j)$ gives the quiver of Example \ref{exmp}.
\end{example}

\subsection{Serre functor}
A Serre functor on $\D$ is an autoequivalence $S_{\D}$ of $\D$ for
which there are natural isomorphisms
\[\Hom_{\D}(A,B)=\Hom_{\D}(B,S_{\D}(A))^*.\]
It is easy to show that if a Serre functor exists then it is unique up to isomorphism.
The motivating example is when $Z$ is a smooth projective variety
of dimension $k$. Then
\[S_{\D}(-)=(-\tensor \omega_Z)[k]\] is a Serre functor on $\D=\D(Z)$.
 A saturated triangulated
category of finite type  always
 has a Serre functor \cite[Corollary 3.5]{bk}.
\smallskip

Suppose $\E\subset \D$ is an exceptional collection, and assume
that $\D$ is saturated. The categories $^{\perp} \E$ and
$\E^\perp$ are then also saturated \cite[Prop. 2.8]{bk} and hence
have Serre functors. Since the mutation functor $\L_\E\colon
{^\perp} \E \lra \E^\perp$ is an equivalence it must commute with
these functors. Thus there is a commutative diagram
\[\begin{CD} {^\perp} \E &@>L_{\E}>> &\E^{\perp} \\
@VS_{{^\perp} \E}VV && @VVS_{\E^\perp}V \\
{^\perp} \E &@>>L_{\E}> &\E^{\perp} \end{CD}\] On the other hand,
the definition of a Serre functor shows that $S_\D$ restricts to
give an equivalence
\begin{equation}
\label{ffs}S_\D|_{{^\perp}\E} \colon {^\perp}\E\lra
\E^\perp.\end{equation} The following Lemma shows that this
functor is the diagonal of the above square.

\begin{lemma}
\label{serre} Suppose  $\E$ is an exceptional
collection in $\D$.  Then
\[S_{\E^\perp}\circ L_\E= S_\D|_{{^\perp}\E} = L_\E \circ S_{{^\perp}\E}.\]
\end{lemma}

\begin{pf}
This is well-known to experts. We give a  proof in Appendix B.
\end{pf}

\begin{cor}
\label{lastly}
If $(E_1,\cdots,E_n)$ is a full exceptional collection in $\D$ then
\[S_\D(E_n)=\L_{E_{n-1}} \cdots \L_{E_{1}}(
E_n).\]
\end{cor}

\begin{pf}
Set $\E=(E_1,\cdots,E_{n-1})$. By Lemma \ref{weedoris} the subcategory ${^\perp}\E\subset \D$ is generated by
a single exceptional object $ E_n$. It follows that it has
trivial Serre functor. Applying Lemma \ref{serre} gives the result.
\end{pf}


\section{Helices}
\label{heli}
 Unlike the last section, the material in this section
is not entirely standard. First we introduce a new and more
flexible definition of a helix which we feel is an improvement on
previous definitions. We then consider rolled-up helix algebras,
and prove Theorem \ref{snow} from the introduction.

\subsection{The definition}
\label{hel} We shall use the following definition of a helix.

\begin{defn}
 A sequence of
objects $\H=(E_i)_{i\in \Z}$ in $\D$ is a \emph{helix} if there exist
positive integers $(n,d)$ with $d\geq 2$ such that
\begin{itemize}
\item[(i)] for each $i\in \Z$ the corresponding \emph{thread}
$(E_{i+1}, \cdots, E_{i+n})$ is a full  exceptional collection,
\smallskip

\item[(ii)] for each $i\in \Z$ one has $E_{i-n}=S_{\D}(E_i)[1-d]$.
\end{itemize}
\end{defn}

Note that by condition (i) we may as well assume that $\D$ has a
full exceptional collection. Then $\D$ is saturated and hence has
a Serre functor $S_{\D}$, so condition (ii) makes sense. The
reason for the apparently strange choice of shift in (ii) will
become clear later.

\begin{remarks}
\label{hlx}
\begin{itemize}
\item[(a)] The pair $(n,d)$ is determined by the helix $\H$;
indeed $n$ is the rank of $K(\D)$ and $d$ is then determined by
condition (ii). We say that $\H$ is of type $(n,d)$.\smallskip

\item[(b)] By Corollary \ref{lastly} condition (ii) is equivalent
to the statement that for all $i\in\Z$
\begin{equation*}
E_{i-n}=\L_{E_{i-(n-1)}} \cdots \L_{E_{i-1}}(
E_i)[1-d].\end{equation*}

\item[(c)] It follows from Remarks \ref{gather}(a) and \ref{hlx}(b)
that it is enough to check that a single thread of $\H$ is a full
exceptional collection.\smallskip

\item[(d)] A full \ec $\E$ of length $n$ generates a helix $\H$ of
type $(n,d)$ for each $d$ in the obvious way. Conversely, a helix
$\H$ is determined by a single thread $\E\subset \H$ together with
the number $d$.
\end{itemize}
\end{remarks}

A helix $\H=(E_i)_{i\in\Z}$ is said to be \emph{geometric} if
for all $i<j$
\[\Hom_{\D}^k(E_i,E_j)=0\text{ unless }k=0.\]
It is \emph{strong} if it satisfies the weaker condition that each thread is a strong exceptional collection.

\begin{example}
\label{eghel} We give examples when $\D=\D(Z)$ with $Z$ a
projective variety.
\begin{itemize}
\item[(a)] Take $Z=\PP^{d-1}$. The canonical bundle is
$\omega_Z=\O(-d)$. The sequence
\[(\cdots ,\O,\O(1),\O(2), \cdots)\]
  is a geometric helix of type $(d,d)$.
\item[(b)] Take $Z=\PP^1\times \PP^1$. The canonical bundle is
$\omega_Z=\O(-2,-2)$. The sequence
\[(\cdots ,\O,\O(1,0),\O(0,1),\O(1,1),\O(2,2),\cdots)\]
is a geometric helix of type $(4,3)$. \item[(c)] Take
$Z=\PP^1\times \PP^1$ again. Then
\[(\cdots ,\O,\O(1,0),\O(3,1),\O(4,1),\O(2,2),\cdots)\]
is a non-strong helix of type $(4,3)$.
\end{itemize}
\end{example}

In previous treatments helices have usually been defined via the
condition in Remark (b), and have
 been required to be of type $(n,n)$.
This means that perfectly sensible helices such as Examples
\ref{eghel}(b),(c) are disallowed.
\smallskip

Later we shall need

\begin{lemma}
\label{du} Suppose $\H$ is a helix of type $(n,d)$ and take a
thread $(E_1,\cdots,E_n)$ with dual collection $(F_n,\cdots,F_1)$.
Then the dual collection to the neighbouring thread $(E_0,\cdots,
E_{n-1})$ is
\[\big(\L_{F_n}(F_{n-1}), \cdots, \L_{F_n}(F_{1}), F_n[1-d]\big).\]
\end{lemma}

\begin{pf}
This is easily proved by directly checking the defining property
\eqref{freddie} of the dual collection.  Indeed, if $1\leq i,j\leq
n-1$ then $\Hom_\D^\blob(E_i,F_n)=0$ so
\[\Hom_\D^\blob(E_i, L_{F_n}(F_j))=\Hom_\D^\blob(E_i,F_j)=\C^{\delta_{ij}}.\]
By Remark \ref{hlx}(b) and \eqref{chee} there is an isomorphism
$E_0=F_n[1-d]$. Since $\L_{F_n}(F_j)\in F_n^{\perp}$ one thus
obtains the remaining vanishing.
\end{pf}

\subsection{Rolled-up helix algebras}

A helix $\H=(E_i)_{i\in\Z}$ of type $(n,d)$ defines a graded algebra
\[A(\H)={\bigoplus_{k\geq 0}}\prod_{j-i=k} \Hom_{\D} (E_i,E_j)\]
known as the \emph{helix algebra}. It has a
$\Z$-action induced by the Serre functor
\[S_\D[1-d]\colon \Hom_{\D}(E_i,E_j) \to \Hom_{\D}(E_{i-n},E_{j-n}).\]
The \emph{rolled-up helix algebra} $B(\H)$ is defined to be the
subalgebra of invariant elements.
\smallskip

Both algebras $A(\H)$ and $B(\H)$ are graded, and it follows from
Lemma \ref{simon} that they are quiver algebras. Given integers
$i<j$ let us write $a_{ij}$ for the dimension of the cokernel of
the map
\[\bigoplus_{i<k<j}\Hom_\D(E_i,E_k)\tensor \Hom_\D(E_k,E_j)\lra\Hom_\D(E_i,E_j).\]
If $i\geq j$ we set $a_{ij}=0$.  The quiver underlying $A(\H)$ has
vertices labelled by the elements of $\Z$ and $a_{ij}$  arrows
connecting vertex $i$ to vertex $j$. The quiver underlying $B(\H)$
has vertices corresponding to elements of $\Z/n\Z$ and \[ n_{ij}=\sum_{p\in\Z} a_{i,j+pn}\]
arrows connecting vertex $i$ to vertex $j$.

Given a thread $\E=(E_{j+1},\cdots,E_{j+n})\subset \H$ there is a bijection between the vertices of the quivers underlying $A(\E)$ and $B(\H)$,
which sends the vertex of $A(\E)$ corresponding to an object $E_i\in\E$ to the vertex of $B(\H)$ labelled by $i\in \Z/n\Z$.
It is then easy to see that the quiver for $B(\H)$ is obtained from that for $A(\E)$
by adding extra arrows corresponding to irreducible morphisms in $\H$ not contained in the thread $\E$.

We will show later (Prop. \ref{plb}) that
in the cases of most interest to us, $a_{ij}=0$ unless $j-i\leq
n$. In that case the extra arrows in $B(\H)$ all point backwards with respect
to the natural ordering on the vertices of $A(\E)$, and, in
particular, the quiver for $B(\H)$ has no loops.

\begin{example}
Take $Z=\PP^1\times \PP^1$. Consider the geometric helix
\[\H=(\cdots ,\O,\O(1,0),\O(0,1),\O(1,1),\O(2,2),\cdots)\]
of Example \ref{eghel}(b). The quiver for the rolled-up helix algebra $B(\H)$ is
easily seen to be
\[
\xymatrix@C=5em{ \bullet\ar[d]_{2}\ar[r]^{2} & \bullet\ar[d]^{2}\\
\bullet\ar[r]_{2} &\bullet \ar[ul]^{4}
 }
\]
Comparing with Example \ref{exmp} we see that the quiver for
$B(\H)$ is obtained from the quiver for the homomorphism algebra
of the thread
\[\E=(\O,\O(1,0),\O(0,1),\O(1,1))\]
by adding extra arrows corresponding to irreducible morphisms in
$\H$ not contained in $\E$.
\end{example}

\subsection{Geometric interpretation}
\label{proof}

Suppose now that $\D=\D(Z)$ with $Z$ a smooth projective Fano
variety. Let $Y$ denote the total space of the canonical line
bundle of $Z$ with its bundle  map $\pi\colon Y\to Z$. Set
$n=\operatorname{rank} K(Z)$ and $d=\dim_{\C}(Y)$. The following
result is adapted from \cite[Prop. 4.1]{br} and \cite[Prop.
7.2]{vdb}.

\begin{thm}
\label{helequiv} Let $B=B(\H)$ be the rolled-up helix algebra of a
geometric helix of type $(n,d)$ on $Z$. Then $B$ is a graded
CY$_d$ quiver algebra which is noetherian and finite over its
centre. Given a thread $\E\subset \H$  there is an equivalence
\[\Phi_{\E}\colon \D(B)\lra \D(Y)\]
sending $B$ to the object $\pi^*(E)$, where  $E=\bigoplus_{E_j\in\E} E_j$.
\end{thm}

%

\begin{pf}
Without loss of generality we can assume that
$\E=(E_1,\cdots,E_n)$. The statement that the collection $\E$ is
full is equivalent to the statement that $E$ classically generates
$\D(Z)$ in the sense of \cite{bvdb}. By \cite[Theorems 2.1.2 and
3.1.1]{bvdb} this is in turn  equivalent to the statement that $E$
generates the category $\DQcoh(Z)$. By the adjunction
\[\Hom^{\blob}_Y(\pi^*(E),F)=\Hom^{\blob}_Z(E,\pi_*(F))\]
and the fact that $\pi_*$ is affine this implies that $\pi^*(E)$
generates $\DQcoh(Y)$. Applying the same argument in reverse this
means that $\pi^*(E)$ classically generates $\D(Y)$.

Next note that $\pi_*(\OO_Y)$ is the sheaf of algebras
\[\A=\pi_*(\OO_Y)=\bigoplus_{p\geq 0} \omega_Z^{-p}\]
with the obvious product structure.  Since $\pi$ is affine it is
then standard (see \cite[Ex. II.5.17]{ha}) that there is an
equivalence of categories between $\OO_Y$-modules on $Y$ and
$\A$-modules on $Z$, defined by sending an $\OO_Y$-module $M$ to
the $\OO_Z$-module $\pi_*(M)$ together with the induced module
structure $\A \tensor \pi_*(M) \to \pi_*(M)$.
 The object $\pi^*(E)$ on $Y$ corresponds under this equivalence  to
\[\pi_* \pi^*(E)= E\tensor_{\OO_Z} \A=\bigoplus_{i\geq 0} E_i\]
with the induced $\A$-module structure.  The assumption that $\H$
is geometric implies that
\[\Ext^p_Y(\pi^*(E),\pi^*(E))=\Ext^p_Z(E,\pi_* \pi^* (E))=0 \text{ for all }p>0.\]
Thus $\pi^*(E)$ is a classical tilting object in the sense of
\cite{hvdb}. The endomorphisms of this object are precisely the
rolled-up helix algebra $B$.

The algebra
\[A=\Gamma(Y,\OO_Y)=\bigoplus_{p\geq 0}\Gamma(Z,\omega_Z^{ -p}).\]
is finitely-generated by the assumption that $\omega_Z^{-1}$ is
ample (see \cite[Theorem 2.3]{reid}).  The affine variety
$X=\operatorname{Spec}(A)$ is the cone over the variety
$Z=\operatorname{Proj}(A)$ and there is a projective morphism
$p\colon Y\to X$ contracting the zero-section of the bundle
$\pi\colon Y\to Z$. Applying \cite[Theorem 7.6]{hvdb} we conclude
that $B$ is finite over its centre and has finite global
dimension, and that there is an equivalence as claimed.

Finally the equivalence $\Phi_\E$ restricts to an
equivalence
\[\D_{\fin}(B) \lra \D_c(Y)\]
where $\D_c(Y)\subset \D(Y)$ is the subcategory of objects with
compact supports. The CY$_d$ property then follows from Serre
duality on $Y$.\end{pf}

Let $Q$ be the quiver underlying the rolled-up helix algebra
$\H=B(\H)$. The vertices are in natural bijection with the objects
of any thread $\E\subset \H$. Let $P_i$ and $S_i$ denote the
projective and simple modules associated to a vertex $i\in Q_0$.
The analogue of Lemma \ref{reallydull} is

\begin{lemma}
\label{andy}
Under the equivalence $\Phi_\E$ of Theorem \ref{helequiv}
one has
\[\Phi_\E (P_i)=\pi^*(E_i), \qquad \Phi_\E(S_i)=i_*(F_i)\]
where  $\F$ be the dual collection to the thread $\E$, and $i\colon Z\into Y$ is the inclusion of the
zero-section.
\end{lemma}

\begin{pf}
The adjunction $\pi^*\dashv \pi_*$ and the fact that $\pi\circ i=\id_Z$ gives
\[\Hom_Y(\pi^*(E_i),i_*(F_j))=\Hom_Z(E_i,F_j).\]
The result then follows as in Lemma \ref{reallydull}.
\end{pf}

\section{Mutation operations}
\label{mutex}

The word mutation
is used to mean two different but related things in the theory of exceptional collections.
Mutation functors were defined in Section 2 and are equivalences
between certain subcategories of our triangulated category $\D$. Mutations of exceptional
collections on the other hand are operations on the set of all
exceptional collections in $\D$. In this section we focus on
mutations in this second sense. The main question that arises is
under what circumstances such mutations preserve purity properties of
exceptional collections.


\subsection{Standard mutations}
 Consider
the  set of all exceptional collections in $\D$ of a certain
length, say $n$. We can use mutation functors to define operations
on this set. The standard way to do this is to define an operation $\sigma_i$ for each
$1<i\leq n$ by the rule
\[\sigma_i(E_1,\cdots,E_{i-2},E_{i-1},E_i,E_{i+1},\cdots, E_n)=\big(E_1,\cdots, E_{i-2}, \L_{E_{i-1}}
(E_i)[-1],E_{i-1}, E_{i+1},\cdots ,E_{n}\big).\]  It is easy to
check  that this operation does indeed take exceptional
collections to exceptional collections. By Remark \ref{gather}(a)
it also takes full collections to full collections.

 \begin{thm}[Bondal, Gorodentsev, Rudakov]
 \label{braid}
 The operations $\sigma_2,\cdots, \sigma_n$ satisfy the braid relations
 \begin{align*}
 \sigma_i\circ\sigma_{i+1}\circ\sigma_i&=\sigma_{i+1}\circ\sigma_i\circ\sigma_{i+1}
 \\
  \sigma_i\circ\sigma_j&=\sigma_j\circ\sigma_i \qquad \text{if }|j-i|>1
  \end{align*}
and hence generate an action of the $n$-string braid group
$\operatorname{Br}_n$ on the set of \ecs of length $n$.
 \end{thm}

\begin{pf}
The second relation is obvious, so we just prove the first. To simplify notation take $n=3$ and $i=2$. Then
\[\sigma_3\circ\sigma_2\circ
\sigma_3(E_1,E_2,E_3)=(L_\F(E_3)[-2],L_{E_2}(E_1)[-1],E_2),\]
where $\F=(E_1,E_2)$. On the other hand
\[\sigma_2 \circ \sigma_3 \circ \sigma_2
(E_1,E_2,E_3)=(L_{\G}(E_3)[-2],L_{E_2}(E_1)[-1],E_2),\] where
$\G=\sigma_2(\F)=(L_{E_2}(E_1)[-1],E_2)$. Applying Remarks
\ref{gather}(a) and (c)  gives the required relation.
\end{pf}

When $\D$ is equipped with a particular choice of t-structure
$\D^{\leq 0}\subset \D$ we shall say that an exceptional
collection $\E$ or a helix $\H$ is \emph{pure} if each constituent
object lies in its heart. The  shift of the mutated object
$\L_{E_{i-1}} (E_i)$ put into the definition of the standard
mutation operation $\sigma_i$ is immaterial for the existence of
the braid group action, since the mutation functors commute with
shifts. However, the given shift ensures that in certain
situations the standard mutations preserve pure
collections. 

\begin{thm}[Bondal, Polishchuk] \label{bop}Suppose $\E=(E_1,\cdots,E_n)$ is a full exceptional collection in $\D$ and take $1<i\leq n$. Suppose
$\D$ is equipped with a t-structure that is preserved by
the autoequivalence $S_{\D}[1-n]$. Then
\begin{itemize}
\item[(i)]  $\E$ pure  implies $\sigma_i(\E)$ pure,\smallskip
\item[(ii)]$\E$ pure  implies $\E$ strong.
\end{itemize}
\end{thm}

Theorem \ref{bop} will be proved in Section \ref{bopo} (as a special case of Theorem \ref{hille}).
It has the following neat consequence.

\begin{cor}
\label{grit}
Suppose $Z$ is a smooth projective variety  satisfying 
\begin{equation*}
\label{bore} \operatorname{rank} K(Z)=\dim(Z) +1.\end{equation*}
Then the standard mutation operations $\sigma_i$ preserve full
collections of sheaves, and all such collections are  strong.
\end{cor}

\begin{pf}
Take $\D=\D(Z)$ equipped with the standard t-structure. The length
of any full exceptional collection in $\D$ is the rank of the
Grothendieck group $K(\D)=K(Z)$.  Since
$S_{\D}(E)=E\tensor\omega_Z[\dim(Z)]$  the result follows  from
Theorem \ref{bop}.
\end{pf}

Corollary \ref{grit} applies for example if $Z=\PP^k$ is a projective space.
On other varieties  the standard mutation operations need
 not preserve collections of sheaves (see Example \ref{cry} below).
 Our aim will be to develop  classes of mutation operations which preserve purity properties in more general situations.


\subsection{Block mutations} Recall that two objects $E$ and $E'$ of $\D$ are said to be \emph{orthogonal} if
\[\Hom_{\D}^\blob(E,E')=0=\Hom_{\D}^\blob(E',E).\]
It is possible to generalise Theorem \ref{bop} by considering
exceptional collections which can be split up into blocks of
mutually orthogonal objects. Such collections were studied by
Hille \cite{hi}.

\begin{defn}
\label{blocks} A \emph{$d$-block} exceptional collection is an exceptional collection
together with a partition of $\E$ into $d$ subcollections
\[\E=(\E_1,\cdots, \E_d),\]
called blocks, such that the objects in each block $\E_i$ are mutually orthogonal.
\end{defn}

For each integer $1<i\leq d$ we can define an operation $\tau_i$ on $d$-block collections in $\D$ by the rule
\[\tau_i \big(\E_1,\cdots,\E_{i-2},\E_{i-1},\E_i,\E_{i+1}\cdots,\E_d\big)=\big(\E_1,\cdots,\E_{i-2},\L_{\E_{i-1}} (\E_i)[-1],
 \E_{i-1}, \E_{i+1},\cdots, \E_d\big).\] Here, if $\E_i=(E_{a+1}, \cdots,
E_b)$ then by definition
\[\L_{\E_{i-1}} (\E_i)=(\L_{\E_{i-1}} (E_{a+1}), \cdots, \L_{\E_{i-1}}(E_b)).\]
 The same proof as before shows that these operations  define an action of the group $\operatorname{Br}_d$ on the set of $d$-block exceptional collections.

\begin{thm} \label{hille}Suppose $\E=(\E_1,\cdots,\E_d)$  is a full $d$-block collection and take $1<i\leq d$.
Suppose
$\D$ is equipped with a t-structure  that is preserved by
the autoequivalence $S_{\D}[1-d]$. Then
\begin{itemize}
\item[(i)] $\E$ pure implies $\tau_i(\E)$ pure,\smallskip
\item[(ii)] $\E$ pure implies $\E$ strong.
\end{itemize}
\end{thm}

We give the proof of Theorem \ref{hille} in Section \ref{bopo}. Note that Theorem \ref{bop} is a special case of Theorem \ref{hille}
since any full exceptional collection of length $n$ is automatically an $n$-block collection, and the operations $\tau_i$ then agree
with the standard mutations $\sigma_i$.

The same argument we gave for Corollary \ref{grit} gives

\begin{cor}
\label{acc} Suppose $Z$ is a smooth projective variety of
dimension $d-1$. Then the $d$-block mutation operations $\tau_i$
preserve full $d$-block collections of sheaves,
 and all such collections are strong.
\qed
\end{cor}

\begin{example}
\label{cry}
 Consider the full, strong exceptional collection of Example \ref{bib}(b) on $Z=\PP^1\times \PP^1$. Then
\[\sigma_3\big(\O,\O(1,0),\O(0,1),\O(1,1)\big)=\big(\O, \O(0,1)[-1],\O(1,0),\O(1,1)\big)\]
which does not consist of sheaves, and
\[\sigma_4\big(\O,\O(1,0),\O(0,1),\O(1,1)\big)=\big(\O, \O(1,0),\O(-1,1),\O(0,1)\big)\]
which consists of sheaves, but is not strong.
If we instead consider the collection as a 3-block collection we have
\[\tau_2\big(\O,[\O(1,0),\O(0,1)],\O(1,1)\big)=\big([\O(-1,0),\O(0,-1)],\O,\O(1,1)\big)\]
which in accordance with Corollary \ref{acc} is strong and consists of sheaves.
\end{example}

\subsection{Mutations of helices}
There are also mutation operations on the set of helices in
$\D$ of a  fixed type $(n,d)$.
 Given such a helix $\H=(E_j)_{j\in \Z}$ and an
integer $i$, we define a new helix
$\sigma_i(\H)=\H'=(E'_j)_{j\in\Z}$ by the rule
\[E'_j=\Biggl\{\begin{array}{lll} E_{j-1} & \text{ if }j= i \mod n, \\ L_{E_j}( E_{j+1})[-1] &\text{ if } j= i-1 \mod n, \\ E_j &\text{ otherwise.}\end{array}\Biggr.\]
Note that $\H'$ satisfies the periodicity property in the
definition of a helix because the Serre functor is an equivalence
and hence commutes with mutations. If $j-i\neq -1 \mod n$ then
\[(E'_{j+1},\cdots, E'_{i-1},E'_i,\cdots, E'_{j+n})=(E_{j+1},\cdots, \L_{E_{i-1}} (E_i)[-1],E_{i-1},\cdots, E_{j+n})\]
Thus  $\H'$ can be obtained by mutating any thread of $\H$ containing
the objects $E_{i-1}$ and $E_{i}$  and then taking the
corresponding helix. The threads of the mutated helix are full
by Remarks \ref{gather}(a) and \ref{hlx}(c).

The periodicity property ensures that the operation $\sigma_i$
only depends on the class of $i$ modulo $n$. There is another
natural operation $\rho$ on helices defined by turning the screw
\[\rho\big[(E_i)_{i\in \Z}\big]=(E_{i+1})_{i\in\Z}.\] Theorem
\ref{braid} immediately implies the following result.

\begin{thm}
\label{braid2} The operations $\rho$ and $\sigma_i$ for
$i\in\Z/n\Z$ satisfy the relations
\begin{align*}
\rho \circ\sigma_i  &= \sigma_{i+1}\circ \rho \\
 \sigma_i\circ\sigma_{i+1}\circ\sigma_i&=\sigma_{i+1}\circ\sigma_i\circ\sigma_{i+1}
 \\
  \sigma_i\circ\sigma_j&=\sigma_j\circ\sigma_i \qquad \text{if }j-i\neq \pm 1
  \end{align*}
 and hence define an action of the affine braid group $\widehat{\operatorname{Br}}_n$ on the set of helices of type $(n,d)$.
\qed\end{thm}

Just as the usual $n$-string braid group $\operatorname{Br}_n$ is
the fundamental group of the configuration space of $n$ points in
$\C$, the affine braid group occurring in Theorem \ref{braid2} is
the fundamental group of the configuration space of $n$ points in
$\C^*$ (see \cite{kp}).

\subsection{Levelled mutations}
Finally in this section we introduce the notion of a levelling,
which will allow us to keep track of operations  defined by mutations
through subcollections.

\begin{defn}
A \emph{levelling} on an \ec $\E=(E_1,\cdots, E_n)$ is a function
$\phi\colon \E \to \Z$ such that $i\leq j \implies \phi(E_i)\leq
\phi(E_j)$.
\end{defn}

A levelling on an \ec $\E$ partitions $\E$ into disjoint
subcollections $\E_i=\phi^{-1}(i)$. We refer to these
subcollections as levels. Note that we do not assume that the
objects in a
given level are orthogonal.\\

There are mutation operations $\sigma_i$ on pairs $(\E,\phi)$ consisting of an exceptional collection $\E$ and a levelling $\phi\colon\E\to\Z$.
Given $i\in\Z$ we define the mutation $\sigma_i(\E,\phi)$ of the pair
$(\E,\phi)$ at level $i$ to be the pair $(\E',\phi')$ for which
the corresponding levels are $\E'_{j}=\E_j$ for $j\notin
\{i-1,i\}$, and
\[\E'_{i-1}=\L_{\E_{i-1}} (\E_i)[-1],\quad \E'_{i} =\E_{i-1}.\]


Suppose for example that $\E=(\E_1,\cdots,\E_d)$ is a $d$-block
collection. There is a unique levelling $\phi\colon \E\to \Z$
satisfying $\phi^{-1}(i)=\E_i$. We call it the \emph{canonical
levelling}. For $1<i\leq d$ the mutated pair
$(\E',\phi')=\sigma_i(\E,\phi)$ is precisely the $d$-block collection
$\tau_i(\E)$ together with its canonical levelling.

\smallskip

We also consider levellings on helices.

\begin{defn}
A levelling on a helix $\H=(E_i)_{i\in\Z}$ of type $(n,d)$ is a function
$\phi\colon \H\to \Z$ such that
\begin{itemize}
\item[(i)] $i\leq j \implies \phi(E_i)\leq \phi(E_j)$,\smallskip
\item[(ii)] for each $i\in \Z$ one has
$\phi(E_{i+n})=\phi(E_i)+d$.
\end{itemize}
\end{defn}


There are mutation  operations $\sigma_i$ defined on pairs
$(\H,\phi)$ consisting of a helix $\H$ of type $(n,d)$ and a
levelling $\phi\colon\H\to\Z$. The levelling  $\phi$  partitions
$\H$ into exceptional collections $\E_j=\phi^{-1}(j)$. To define
$\sigma_i (\H,\phi)$ take a thread $\E$  of $\H$ containing the
levels $\E_i$ and $\E_{i-1}$. This can always be done by the
periodicity of $\phi$ and our assumption that $d\geq 2$ in the
definition of a helix. Since $\phi$ restricts to a levelling on
$\E$, we can define the mutated pair
\[(\E',\phi')=\sigma_i(\E,\phi)\] and then set $\sigma_i(\H,\phi)$
to be the helix generated by $\E'$, together with the unique
levelling extending $\phi'$. It is easy to check that the
resulting mutated helix $\sigma_i(\H,\phi)$ is independent of the
thread $\E$ we chose, and depends only on
 the value of $i$ modulo $d$.


\section{An argument of Bondal and Polishchuk}

\label{bopo} In this section we give a proof of Theorem
\ref{hille}. This involves extending several arguments from
\cite{bo,bp} to the $d$-block case. We will later use a similar
argument to prove our main result Theorem \ref{big}. The crucial
ingredient in the proofs of these theorems is the following simple
observation from \cite{bo}.

\begin{lemma}
\label{round} Let $(\H,\phi)$ be a pair of a helix $\H$ of type
$(n,d)$ together with a levelling $\phi\colon \H \to \Z$. Suppose
also that the objects of $\E_m=\phi^{-1}(m)$ are mutually
orthogonal. Then
\[(\sigma_{m-(d-2)}\circ \cdots \circ \sigma_{m})(\H,\phi)=(\H,\phi+1).\]
\end{lemma}

\begin{pf}
Adding a constant to the levelling we may as well assume that
$m=d$. Consider the thread $\E=(\E_1, \cdots, \E_{d-1}, \E_d)$ and
also the subcollection $\E'= (\E_{1}, \cdots, \E_{d-1})$. Applying
the mutations in the statement of the Lemma to the thread $\E$ we
obtain the collection
\[(\L_{\E'}(\E_d)[1-d],\E_{1}, \cdots, \E_{d-1}).\]
Lemma \ref{weedoris} implies that ${^\perp}\E'=\langle\E_d\rangle$
and this category has trivial Serre functor by the assumption that
all objects of the collection $\E_d$ are orthogonal. Thus by Lemma
\ref{serre} and the periodicity of the helix
\[\L_{\E'}(\E_d)[1-d]=S_\D(\E_{d})[1-d]=\E_{0}.\]
The mutated thread is therefore \[(\E_{0}, \E_{1}, \cdots ,
\E_{d-1})\] and so the helix it generates is $\H$. The induced
levelling is easily seen to be $\phi+1$.
\end{pf}

A simple consequence is

\begin{lemma}
\label{roundtwo} Let $(\H,\phi)$ be a pair of a helix $\H$ of type
$(n,d)$ together with a levelling $\phi\colon \H \to \Z$. Suppose
also that the objects of $\E_d=\phi^{-1}(d)$ are mutually
orthogonal. Then for any integer $1\leq k\leq d$
 \[\L_{\E_{k}} \cdots \L_{\E_{d-1}} (\E_d) [k-d]=
\R_{\E_{k-1}}\cdots \R_{\E_{1}}(\E_0)[k-1]\] \end{lemma}

 \begin{pf}
Note that Lemma \ref{round} with $m=d$ implies
\begin{equation}
\label{who}(\sigma_{k+1}\circ \cdots \circ
\sigma_{d})(\H,\phi)=(\sigma_2\circ \cdots
\sigma_k)^{-1}(\H,\phi-1).\end{equation} Now the inverses of the
mutation operations $\sigma_i$ can be written in terms of right
mutations in the obvious way. Taking the $k$th level of both sides
of \eqref{who} then gives the result.
\end{pf}

Define
a block structure on a helix $\H$ to be a levelling $\phi\colon \H
\to \Z$ such that any two objects in the same level $\phi^{-1}(m)$
are orthogonal.
Clearly, the $(n,d)$ helix $\H$ generated by a $d$-block
exceptional collection $\E$ has a block structure defined by
extending the canonical levelling from $\E$ to $\H$.
 Conversely, a block structure
on a helix of type $(n,d)$ partitions each thread of $\H$ into at most $d+1$
blocks.

\begin{thm}
\label{di}
Let $\D^{\leq 0}\subset\D$ be a t-structure and
suppose $\H$ is a helix in $\D$ with a block structure $\phi\colon\H\to\Z$. Take an integer $m\in\Z$ and set $(\H',\phi')=\sigma_m(\H,\phi)$.
Then
\begin{itemize}
\item[(i)] $\H$ pure implies $\H'$ pure,\smallskip \item[(ii)]
$\H$ pure implies $\H$ geometric,\smallskip \item[(iii)] $\H$
geometric implies $\H'$ geometric.
\end{itemize}
\end{thm}

\begin{pf}
Adding a constant to the levelling we may as well assume that $m=d$. Suppose first that $\H$ is pure and consider the thread
\[\E=(\E_{1},\cdots,\E_{d-1},\E_d).\]
By periodicity, to prove (i) it will be enough to show that the
mutated thread
\[\tau_d(\E)=(\E_{1},\cdots,\L_{\E_{d-1}} (\E_d) [-1],\E_{d-1})\]
is pure.
In fact we prove more, namely that for $0\leq k<d$ the multiply mutated thread
\[\E'=\tau_k \cdots \tau_d(\E)=(\E_{1},\cdots,\E_{k-1}, \L_{\E_{k}} \cdots \L_{\E_{d-1}} (\E_d) [k-d],\E_{k},\cdots \E_{d-1})\]
is pure.

Consider the mutation functor $\L_{\E_j}$ applied to an element
$X\in{^\perp}\E_j$. The assumption that all objects of $\E_{j}$
are orthogonal implies that there is a triangle
\begin{equation}
\label{cold}\L_{\E_{j}} (X)[-1] \lra \bigoplus_{E\in\E_{j}}
\Hom_{\D}^\blob(E,X)\tensor E \lra X.\end{equation} This can be
seen by writing the functor $\L_{\E_{j}}$  as a repeated mutation
and using the octahedral axiom. Suppose now that $X\in\D^{\geq
0}$. Since $\E$ is pure, the graded vector space
$\Hom_{\D}^\blob(E,X)$ is concentrated in  non-negative degrees,
and so using the triangle \eqref{cold}
 we can conclude that $\L_{\E_j}(X)[-1]\in\D^{\geq 0}$ also.

Applying this result repeatedly shows that $\E'\subset \D^{\geq
0}$. Now note that by Lemma \ref{roundtwo}  the thread $\E'$ can
also be written as a repeated right mutation. A similar argument
to the above then gives $\E'\subset \D^{\leq 0}$ and hence $\E'$
is pure as claimed.

To prove (ii) assume $\H$ is pure and take a thread $\E\subset\H$
that is partitioned by $\phi$ into $d$ blocks.  We claim that $\E$
is strong. It then follows from  Proposition \ref{below} below
that $\H$ is geometric. To prove the claim suppose that there are
elements $E_i$ and $E_j$ of $\E$  with $\Hom_{\D}^k(E_i,E_j)\neq
0$ for some $k>0$. By applying block mutations we can move blocks
between $E_i$ and $E_j$ out of the way and arrive at a  collection
such that $E_i$ and $E_j$ lie in neighbouring blocks $\E_{m-1}$
and $\E_m$ say. Applying one more mutation gives an exceptional
collection with blocks
\[(L_{\E_{m-1}}(\E_m)[-1],\E_{m-1}).\]
By part (i) this collection is pure. But by the triangle
\eqref{cold} above
\[\Hom_\D^{-k}(L_{\E_{m-1}}(E_j)[-1],E_i)=\Hom_\D^k(E_i,E_j)^*\neq
0\] giving a contradiction.

To prove (iii) let $\D^{\leq 0}\subset \D$ be the t-structure
corresponding to the standard t-structure on $\D(A)$ under the
equivalence $\Phi_\E$ of Theorem \ref{RB}.
 By Remark \ref{extra} the
 assumption that $\H$ is geometric implies that every object far enough to the right in $\H$ is pure.
 The argument for part (i)  shows that the same is true for the mutated helix $\H'$.
 The argument we gave for (ii) then shows that $\H'$ is geometric.
\end{pf}

\begin{prop}
\label{below}
Suppose $\E=(\E_1,\cdots,\E_d)$ is a full, strong $d$-block collection.
Let $\H$ be the helix of type $(n,d)$ generated by $\E$. Then $\H$ is geometric iff
\[\Hom_{\D}^k(E_i,E_j)=0\text{ for all }i<j\text{ and }k<0.\]
\end{prop}

\begin{pf}
Let $A=A(\E)$ be the homomorphism algebra of $\E$. Applying Remark
\ref{extra} to the opposite category $\D^{\op}$ we see that there
is a t-structure $\D^{\leq 0}\subset \D$ whose heart is a finite
length abelian category $\A$ consisting of objects $X$ for which
$\Hom_\D(X,E_i)$ is concentrated in degree 0 for all $i$. The
simple objects $(G_n,\cdots,G_1)$ of $\A$ satisfy
\begin{equation}
\Hom^{k}_{\D}(G_j,E_i)=\left\{
\begin{array}{ll}
\C &\text{$i=j$ and $k=0$}, \\
0 &\text{otherwise.}
\end{array}\right.\end{equation}
It follows from the definition of the Serre functor that
$F_j=S_\D(G_j)$ where $(F_n,\cdots,F_1)$ is the dual collection to
$\E$.

Let $\tau$ denote the functor $S_\D[1-d]$. We first claim that
$\tau$ is left exact with respect to this t-structure, that is
that $\tau(\D^{\geq 0})\subset \D^{\geq 0}$. To prove this it is
enough to show that $\tau(G_j)=F_j[1-d]\in \D^{\geq 0}$ for each
$j$.

By Theorem \ref{di}(i) the collection
\[\big(\L_{\E_1} \cdots \L_{\E_{d-1}}(\E_d)[1-d], \cdots, \L_{\E_{1}}(\E_2)[-1],\E_1\big)\]
is pure. Note that if $X$ and $Y$ are orthogonal  $\L_X(Y)=Y$. By
\eqref{chee} it follows that the dual collection
$(F_n,\cdots,F_1)$ is a reordering of the collection
\[\big(\L_{\E_1} \cdots \L_{\E_{d-1}}(\E_d), \cdots, \L_{\E_{1}}(\E_2),\E_1\big).\]
Hence $F_j[1-d]$ lies in $\D^{\geq 0}$ for all $j$.

Now suppose the condition of vanishing of negative $\Hom$ groups
in $\H$ holds. By periodicity of $\H$, to prove that $\H$ is
geometric it will be enough to show that for $E,E'\in\E$ and all
$m\geq 0$ one has
\[\Hom_\D^k(\tau^m(E'),E)=0 \text{ for }k>0.\]
Since $E$ is an injective object in the heart of $\D^{\leq 0}$ and
$\tau^m(E')\in \D^{\geq 0}$ the result follows.
\end{pf}

\smallskip

\noindent {Proof of Theorem \ref{hille}.}
 Since the t-structure is preserved by the functor $S_\D[1-d]$,
the helix $\H$ of type $(n,d)$ generated by $\E$ is pure. The
canonical levelling on the $d$-block collection $\E$ extends to a
block structure $\phi\colon\H\to\Z$. Thus (i) follows from Theorem
\ref{di}(i). Part (ii) follows as in the proof of Theorem
\ref{di}(ii). \qed


\section{Height functions}

In this section we introduce special types of levellings called height functions. In the next section we will relate mutations defined by such levellings to the tilting operation for homomorphism and rolled-up helix algebras.

\subsection{Relatedness}

It will be convenient to introduce the following terminology.

\begin{defn}
Let $\E=(E_1,\cdots,E_n)$ be a full exceptional collection with
dual collection $\F=(F_n,\cdots,F_1)$. Let $p\geq 0$. We say that
objects $E_i,E_j\in \E$ with $i\leq j$ are $p$-related in $\E$ if
\[\Hom^k_{\D}(F_j,F_i)=0\text{ for }k\neq p.\]
Note that if $F_j$ and $F_i$ are orthogonal then $E_i$ and $E_j$ are $p$-related for all $p$.
\end{defn}

A simple Corollary of Theorem \ref{hille} is

\begin{lemma}
\label{peewee} Suppose $\E=(\E_1,\cdots,\E_d)$ is a full $d$-block
collection which is pure in some t-structure preserved by
$S_\D[1-d]$. Let $\phi\colon \E \to \Z$ be the canonical
levelling. Then for any $i\leq j$ the objects $E_i$ and $E_j$ are
$(\phi(j)-\phi(i))$-related.
\end{lemma}

\begin{pf}
By Theorem \ref{hille} the collection $\E$ is strong. As we observed in the proof of Proposition \ref{below}  the collection
\[\E'=(\L_{\E_1}\cdots \L_{\E_{d-1}}(\E_d)[1-d],\cdots, \L_{\E_1}(\E_2)[-1],\E_1)\]
is strong, and differs from the dual collection $\F$ by reordering
and the shifts. The result follows.
\end{pf}

Later we shall need

\begin{lemma}
\label{threads} Let $\H=(E_i)_{i\in\Z}$ be a helix of type $(n,d)$.
\begin{itemize}
\item[(a)]  Suppose $i\leq j$ and take two threads
$\E_0,\E_1\subset \H$ containing the objects $E_i$ and $E_j$. Then
$E_i$ and $E_j$ are $p$-related in the thread $\E_0$ precisely if
the same is true in $\E_1$.
\smallskip \item[(b)] Consider the threads
 $\E_0= (E_0,\cdots, E_{n-1})$ and $\E_1=(E_1,\cdots, E_{n})$ and take an integer $1\leq i\leq n-1$.
 Then $E_0$ and $E_i$
 are $p$-related in  $\E_0$ precisely  if $E_i$ and $E_n$ are  $(d-p)$-
 related in $\E_1$.
 \end{itemize}
\end{lemma}

\begin{pf}
For (a) it is enough to consider the case when \[\E_0=(E_0,\cdots,
E_{n-1}),\quad\E_1=(E_1,\cdots, E_{n}),\] and $1\leq i\leq j\leq n-1$.
 Write $\F_1=(F_n, \cdots, F_1)$ for the dual collection to $\E_1$. By Lemma \ref{du} the dual
collection to $\E_0$ is
\[\F_0=(\L_{F_n}(F_{n-1}), \cdots, \L_{F_n}(F_{1}), F_n[1-d]).\]
Now the dual objects to $E_i$ and $E_j$ in $\F_0$ and in $\F_1$
differ by mutation by $F_n$, and since the  mutation functor defines an
equivalence ${^\perp} F_n\to F_n^{\perp}$, the notion of $p$-related
is the same in each case.

For (b) note that  as elements of $\E_1$,  the dual of $E_i$ is
$F_i$, and the dual of $E_n$ is $F_n$.  On the other hand, as
elements of $\E_0$ the dual of $E_0$ is $E_0=F_n[1-d]$ and the
dual of $E_i$ is $\L_{F_n}(F_{i})$. Applying the functor
$\Hom_\D(-,F_n)$ to the triangle \eqref{left} shows that
\[\Hom_{\D}^{-k}(F_n,F_{i})^*=\Hom_{\D}^{k+1}(\L_{F_n}(F_{i}),F_n)\]
and the result follows.
\end{pf}

\subsection{Height functions}
We now make the following crucial definition.

\begin{defn}
\label{height} Let $\E=(E_1,\cdots,E_n)$ be a full, strong
exceptional collection. A levelling $\phi\colon \E \to \Z$ is said
to be \emph{tilting at level $m$} if the following condition
holds. Suppose  $E_i,E_j\in \E$ with $\phi(E_i)=m$ and set
$\phi(E_j)=p$. Then
\begin{itemize}
\item[(i)]  $i\leq j$ implies that $E_i$ and $E_j$ are
$(p-m)$-related, \item[(ii)] $j\leq i$ implies that $E_j$ and
$E_i$ are $(m-p)$-related.
\end{itemize}
A \emph{height function} for an object $E\in \E$
is a levelling $\phi\colon \E \to \Z$ which is tilting  at level
$0$ and which satisfies $\phi^{-1}(0)=\{E\}$.
\end{defn}

In general a height function for an object $E$ in an exceptional
collection $\E$  may or may not exist, and  will usually not be in
any sense unique.

\begin{examples}
\begin{itemize}
\item[(a)] If $Z$ is a smooth projective variety of dimension
$d-1$, and  $\E$ is a full $d$-block collection of sheaves
on $Z$ then by Lemma \ref{peewee} the canonical levelling on $\E$
is tilting at all levels. It follows from Lemma \ref{why} that
height functions exist for all objects of $\E$.

 \item[(b)]
 Consider the full, strong collection
\[\E=(\O,\O(1,0),\O(1,1),\O(2,1))\]
on $Z=\PP^1\times \PP^1$ with dual collection
\[\F=(\O(0,-1)[2], \O(1,-1)[1], \O(-1,0)[1],\O).\]
 The possible
height functions for $\O(2,1)$ take values $(-2,b, -1,0)$ for
$b\in\{-2,-1\}$. Similarly, the possible height functions for
$\O(1,0)$ take values $(-1,0,1,a)$ with  $a\geq 1$.
\smallskip
\item[(c)] In Section \ref{delpezzo} we shall adapt an argument of Herzog to show that if $\E$ is
a strong exceptional collection on a del Pezzo surface,
 then it is possible to reorder $\E$ so that height functions exist for all objects $E\in \E$.

\end{itemize}
\end{examples}

\begin{lemma}
\label{sunn} Suppose $\E=(E_1,\cdots,E_n)$ is a full strong
exceptional collection with dual collection $\F=(F_n,\cdots,F_1)$.
Suppose $\phi\colon\E\to\Z$ is a levelling which is tilting at
level $m$, and take distinct objects $E_i,\E_j\in\E_m$. Then $E_i$
is orthogonal to $E_j$, and moreover $F_i$ is orthogonal to $F_j$.
\end{lemma}

\begin{pf}
By definition of tilting at level $m$ one has
\[\Hom^k_{\D}(F_i,F_j)=0=\Hom^k_\D(F_j,F_i)\]
unless $k=0$. But since the collection $\E$ is strong  the objects $F_i$ and $F_j$
are identified under the equivalence $\Phi_\E$ with simple modules
for the homomorphism algebra $A(\E)$. It follows that they must be orthogonal.

Consider the homomorphism algebra $A(\E)$ and the corresponding
quiver $Q$. Let the vertex of $Q$ corresponding to an object $E_p$
be denoted $p$. By \eqref{mortlock} the number of arrows $n_{pq}$
from vertex $p$ to vertex $q$ is the dimension of the space
$\Hom^1_\D(F_q,F_p)$. In particular if $n_{pq}>0$ then $p<q$. On
the other hand the space $\Hom_\D(E_p,E_q)$ is spanned by paths
in the quiver from vertex $p$ to vertex $q$, modulo relations. By
the first part there are no arrows between vertices corresponding
to elements of $\E_m$. Thus there are no paths between the
corresponding vertices either.
\end{pf}

\begin{cor}
\label{why} Suppose $\E=(E_1,\cdots,E_n)$ is a full, strong
exceptional collection and take an integer $1\leq i\leq n$.
Suppose there is a levelling $\phi\colon \E\to\Z$ which is tilting
at level $m=\phi(E_i)$. Then there is a height function for
$E_i\in\E$.
\end{cor}

\begin{pf}
 Replacing $\phi$ by $\phi-m$ we can assume that $m=0$.
 If $E_j\in \E_0$ with $j<i$ then by Lemma \ref{sunn}
 $E_j$ is $p$-related to $E_i$ for all $p$. In this case we redefine $\phi(E_j)$ to be $-1$.
 Similarly for $E_j\in \E_0$ with $j>i$ we set $\phi(E_j)=+1$. The resulting levelling
 is a height function for $E_i\in \E$.
\end{pf}

%

\subsection{Height functions on helices}

In this section we  extend the definition of height functions to
helices.

\begin{lemma}
Let $\H$ be a strong helix with a levelling $\phi\colon \H \to
\Z$. Suppose $\phi$ restricted to some thread $\E\subset \H$
containing the exceptional collection $\E_m=\phi^{-1}(m)\subset
\H$ is tilting at level $m$. Then $\phi$ restricted to any such
thread is tilting at level $m$.
\end{lemma}

\begin{pf}
We can clearly reduce to the case $m=0$. Take two threads
$\E,\E'\subset \H$ containing $\E_0$. We must show that if
$\phi|_{\E}$ is tilting at level $0$ then so is  $\phi|_{\E'}$. It
is enough to consider the case when $\E$ and $\E'$ are
neighbouring threads, so without loss of generality  we can take
\[\E= (E_0,\cdots, E_{n-1}), \quad \E'=(E_1,\cdots, E_{n}).\]
Assume $\phi|_{\E}$ is tilting at level 0. Take two objects
$E_i,E_j\in \E'$ with $\phi(E_i)=0$ and set $q=\phi(E_j)$. If
$j\leq i$ then $q\leq 0$ and we must show that $E_j$ and
 $E_i$ are $|q|$-related in $\E'$. But this follows from Lemma \ref{threads}(a) because both
$E_i$ and $E_j$ lie in $\E$ and  $\phi|_{\E}$ is tilting at level
0. So suppose $j\geq i$. Then $q\geq 0$ and we must show that
$E_i$ and $E_j$ are $q$-related. If $E_j\in\E$ the same argument
applies. The only other possibility is that $j=n$. Since
$\phi|_{\E}$ is tilting at level 0 it follows that $E_0$ and $E_i$
are $(d-q)$-related.  Applying Lemma \ref{threads}(b) we conclude
that $E_i$ and $E_n$ are $q$-related as required.
\end{pf}

Now we can make the following definition.

\begin{defn}
Let $\H$ be a strong helix.
 A levelling $\phi\colon \H \to\Z$ is \emph{tilting at level $m$} if there is a thread $\E\subset\H$
  containing $\E_m=\phi^{-1}(m)$ such that the restriction of
 $\phi$ to $\E$ is tilting at level $m$.

A \emph{height function} for an object $E\in \H$ is a levelling $\phi\colon \H\to \Z$ which is
  tilting at level $0$ and satisfies $\phi^{-1}(0)=\{E\}$.
\end{defn}

The following result shows that providing one knows that height functions exist on strong exceptional collections
in a given category $\D$ then one can always construct them on strong helices.

\begin{lemma}
\label{pf}
 Let $\H$ be a strong helix and take an object $E\in \H$.  Let
$\E\subset \H$ be the thread which has $E$ at its start. Suppose
there is a height function $\phi\colon \E \to \Z$ for $E\in \E$.
Then  there is a height function for $E\in \H$.
\end{lemma}

\begin{pf}
Label the helix so that $E=E_1$ and therefore
$\E=(E_1,\cdots,E_n)$. By hypothesis there exist height functions
$\phi\colon \E \to \Z$ for $E\in \E$. Choose one such that
$\phi(E_n)$ takes the minimal possible value, say $p$. We claim
that $p\leq d-1$. Extending $\phi$ periodically then gives a height
function for $E\in \H$.

Let $\F=(F_n,\cdots,F_1)$ be the dual collection to $\E$. Note
that $F_1=E_1$. Take $1\leq m\leq n$ such that
$\Hom_\D(E_m,E_{n+1})\neq 0$ and $\Hom_\D(E_j,E_{n+1})=0$ for
$m<j\leq n$. This is possible because $\E$ is full and $\H$ is
strong.

Take $1\leq j\leq n$. By \eqref{chee} and Remark \ref{hlx}(b) the
object $E_{j-n}[d-1]$ is obtained by mutating the object
$F_{j}=\L_{E_1}\cdots \L_{E_{j-1}}(E_j)$ through the collection
$(E_{j+1-n},\cdots,E_{0})$. If $m\leq j\leq n$ these objects have
no maps to $E_1$ by periodicity of $\H$, so
\[\Hom_\D^\blob(E_j,E_{n+1})[1-d]=\Hom_\D^\blob(E_{j-n}[d-1],E_1)=\Hom_\D^\blob(F_{j},F_1).\]
For $m<j$ this implies that $F_{j}$ is orthogonal to $F_1$ and the
value of $\phi(E_j)$ is unconstrained. Since $\phi$ is assumed
minimal this implies that $\phi(E_j)=\phi(E_m)$ for $j>m$. By
assumption $\Hom_\D(E_m,E_{n+1})$ is nonzero, so
$\Hom^\blob_\D(F_m,F_1)$ is nonzero in degree $d-1$. Hence $E_1$
and $E_m$ are $d-1$ related, and $\phi(E_m)=d-1$.
\end{pf}


\section{Tilting and mutation}

This section contains our main results. We relate tilting to mutations defined by height functions. We also show that
mutations determined by height functions preserve geometric helices.

\subsection{Tilting homomorphism algebras}


Let $\E=(E_1,\cdots,E_n)$ be a full, strong exceptional collection
with  homomorphism algebra $A=A(\E)$. Let $Q$ be the underlying
quiver of $A(\E)$,  let $i\in Q_0$  denote  the vertex corresponding to an object
$E_j$ and let $S_j$ be the corresponding simple
$A$-module.

Fix a vertex $i\in Q_0$. By Lemma \ref{min} there is a projective
resolution of the module $S_i$ of the form
\[0\lra \bigoplus_{j\in Q_0} P_j ^{\oplus d_{j,i
}^n} \lra \cdots  \lra \bigoplus_{j\in Q_0} P_j^{\oplus d_{j,i}^1}
\lra P_i \lra S_i\lra 0\] where $d_{j,i}^p=\dim_\C
\Ext^p_{A}(S_i,S_j).$ Given an integer $k\geq 0$ consider the
truncated complex \begin{equation} \label{balls} 0\lra
\bigoplus_{j\in Q_0} P_j ^{\oplus d_{j,i }^k} \lra \cdots  \lra
\bigoplus_{j\in Q_0} P_j^{\oplus d_{j,i}^1} \lra P_i\lra
0\end{equation} situated in degrees $-k$ up to $0$. We define
$P_i^{(k)}$ to be the corresponding element of $\D(A)$.

\begin{prop}
\label{han} Suppose $\phi\colon \E\to \Z$ is tilting at level $m$
and  $E_i\in \E_m$. Then under the equivalence $\Phi_\E$ of
Theorem \ref{RB} one has
\[\Phi_\E(P_i^{(k)})=L_{\E_{m-k}} \cdots L_{\E_{m-1}} (E_i)\]
for all $k\geq 0$.
\end{prop}

\begin{pf}
Let $\F=(F_n,\cdots,F_1)$ be the dual collection to
$\E=(E_1,\cdots,E_n)$. Recall from Lemma \ref{reallydull} that
under the equivalence $\Phi_\E$ the projective module $P_i$
corresponds to $E_i$ and the simple module $S_i$ corresponds to
$F_i$. Since $\phi$ is tilting at level $m$ we have $d_{j,i}^p=0$
unless $E_j\in \E_{m-p}$.

Let $\E'$ be the subcollection
\[\E'=(\E_{m-k},\cdots,\E_{m-1})\subset \E\]
 and set $R=\Phi_\E(P_i^{(k)})\in \D$. We must
show that $R=L_{\E'}(E_i)$.

By definition of the complex $P_i^{(k)}$ there is a triangle \[E_i
\lra R \lra A\] with $A$ in the subcategory $\langle\E'\rangle
\subset \D$.  Since the complex $P_i^{(k)}$ is a truncation of a
projective resolution of the simple module $S_i$ there is also a
triangle
\[R \lra F_i \lra B\]
with  $B$ in the subcategory generated by the elements $E_j\in \E$
with $\phi(E_j)<m-k$. This implies that $R\in \E'^{\perp}$. The
result then follows from Remark \ref{gather}(b).\end{pf}

\begin{cor}
\label{blow}
 Suppose $\phi\colon \E\to \Z$ is tilting at level $m$
and $E_i\in \E_m$. Then the object \[L_{\E_{m-k}} \cdots
L_{\E_{m-1}} (E_i)[-k]\] lies in the extension-closed subcategory
of $\D$ generated by objects of the form $E[-j]$ for $j\geq 0$ and
$E\in (\E_{m-k},\cdots,\E_{m})$. Similarly, the object
\[\R_{\E_{m+k}}\cdots R_{\E_{m+1}}(E_i)[k]\]
lies in the extension-closed subcategory of $\D$ generated by
objects of the form $E[j]$ with $j\geq 0$ and $E\in
(\E_{m},\cdots,\E_{m+k})$.
\end{cor}

\begin{pf}
The first statement is obvious from the proof of Proposition
\ref{han} and the explicit form of the complex \eqref{balls}. The
second statement is most easily obtained by considering the
opposite category $\D^\op$.
\end{pf}

When $k=1$ the object $P_i^{(1)}[-1]$ is precisely the
complex $R_i$ from Section \ref{saville}. It follows that mutations of exceptional collections defined by height functions
induce vertex tilts of the corresponding homomorphism algebras.

\begin{prop}
\label{blbl}
Suppose $\phi\colon \E\to \Z$ a
 height function for $E_i\in \E$.  Set
$(\E',\phi')=\sigma_0(\E,\phi)$ and assume the collection
$\E'$ is strong. Then the homomorphism algebra $A(\E')$ is the left tilt of the algebra $A(\E)$
at the vertex $i$.
\end{prop}

\begin{pf}
Let $Q$ and $Q'$ be the
quivers underlying $A=A(\E)$ and $A'=A(\E')$ respectively. The
vertices of these quivers correspond to the elements of the
collections $\E$ and $\E'$. Set $E'=L_{\E_{-1}}(E)[-1]$. Since
$\E\setminus \{E\}=\E'\setminus\{ E'\}$  there is an obvious
choice of bijection $\psi\colon Q_0\to Q'_0$. Consider the diagram
of equivalences
\[\begin{CD} \D(A) &@>\Psi>> &\D(A') \\
@V\Phi_{\E}VV && @VV\Phi_{\E'}V \\ \D &@>\id>> &\D
\end{CD}\]
Let $i\in Q_0$ be the vertex corresponding to the object $E$.
Applying Lemma \ref{reallydull} we see that if $j\in Q_0$ and
$j\neq i$ then  $\Psi(P_j)=P'_{\psi(j)}$. On the other hand
Proposition \ref{han}  shows that $\Phi_\E(R_i)=E'$ and hence
$\Psi(R_i)=P_{\psi(i)}$.
\end{pf}

\subsection{Preservation of geometric helices}

Our main result shows that geometric helices are preserved under mutations defined by height functions.

%
%
%

\begin{thm}
\label{big} Suppose $\H$ is a geometric helix in $\D$ and $\phi
\colon \H \to \Z$ is a levelling that is tilting at level $m$. Set
$(\H',\phi')=\sigma_m (\H,\phi)$.  Then $\H'$ is geometric.
 \end{thm}

\begin{pf} Adding a constant to the levelling we may as well assume that
$m=d$. Take a t-structure $\D^{\leq 0}\subset \D$ in which all
objects $E\in\H$ satisfying $\phi(E)\geq 0$ are pure. By Remark \ref{extra}
and the fact that $\H$ is geometric, one can construct such a
t-structure by considering the equivalence of Theorem \ref{RB}
corresponding to a thread of $\H$ far enough to the left.

Consider the thread $\E=(\E_1,\cdots, \E_d)\subset \H$ and the mutated thread
\[\E'=(\E_1,\cdots,\E_{d-2},\L_{\E_{d-1}}(\E_{d})[-1],\E_{d-1})\subset \H'.\]
By Lemma \ref{roundtwo} one has
\[\G=\L_{\E_{d-1}}(\E_{d})[-1]=\R_{\E_{d-2}}\cdots \R_{\E_{1}}(\E_0)[d-2].\]
The two parts of
Corollary \ref{blow} then show that $\G\subset \D^{\geq 0}$ and
$\G\subset \D^{\leq 0}$. Thus $\G$ is pure.
 By periodicity of $\H$ it follows that
 \[\Hom_{\D}^{<0}(E_i,E_j)=0\text{ for all }i<j.\]

Take an object
$A\in \G$. The first part of Corollary \ref{blow} implies that $A$ is in the extension-closed subcategory generated by
right shifts of objects $E\in\H$ with $\phi(E)\leq d$. Since $\H$ is geometric it follows that
 \begin{equation}
 \label{mi}
 \Hom_{\D}^{>0}(A,E_j)=0 \text{ if }\phi(E_j)> 0.\end{equation}
Viewing $\G$ instead as a right mutation, the second part of
Corollary \ref{blow} shows that $A$ is in the extension-closed subcategory generated by
left shifts of objects $E\in\H$ with $\phi(E)\geq 0$. This gives \begin{equation}\label{ni}\Hom_{\D}^{>0}(E_j,A)=0\text{ if
}\phi(E_j)< d.\end{equation} Finally put $\kappa=S_\D[1-d]$ and take $B\in
\kappa^r(\G)$ for some $r>0$. By periodicity of the helix,  $B$ is in the extension-closed
subcategory of $\D$ generated by right shifts of objects $E\in\H$
with $\phi(E)\leq 0$. It follows that $\Hom_\D^{>0}(B,A)=0$. This concludes
the proof that the helix generated by $\E'$ is geometric.
\end{pf}

\subsection{Tilting rolled-up helix algebras} In this section we show that mutations of geometric helices defined by
height functions induce vertex tilts of the corresponding  rolled-up helix algebras. First we prove

\begin{prop}
\label{plb}
Let $\H$ be a geometric helix of type $(n,d)$ with $d\geq 2$ and
let $B(\H)$ be the rolled-up helix algebra. Suppose $\phi\colon\H\to\Z$ is a height function for $E_i\in\H$.
Then the quiver $Q$ underlying $B(\H)$ has no loops at the vertex corresponding to $E_i$.
If $d\geq 3$ then there are no 2-cycles in $Q$ passing through this vertex either.
\end{prop}

\begin{pf}
Given an integer $k<i$ we claim that every morphism $E_k\to E_i$ factors via
an element of $\E_{-1}$.
To see this consider the triangle
\begin{equation}
\L_{\E_{-1}} E_i[-1]
\lra \bigoplus_{E\in\E_{-1}} \Hom_{\D}^\blob(E,E_i)\tensor E \lra E_i.\end{equation}
The existence of such a triangle follows as in the proof of Theorem \ref{di}, or from Proposition \ref{han} together with the explicit form of the complex
\eqref{balls}.

If there
is a nonzero morphism $E_k\to E_i$ which does not factor through
an element of $\E_{-1}$ then it gives rise to a nonzero  element of $\Hom^1_\D(E_k,\L_{\E_{-1}}(E_i)[-1])$.
Since $\phi(E_k)<0$ this contradicts equation \eqref{ni} in the proof of Theorem \ref{big}.
Considering right mutations in a similar way shows that if $k>i$ every morphism $E_i\to E_k$  factors via an element of $\E_1$.

Consider now the quiver underlying $B(\H)$. Its vertices can be
put in bijection with the elements of the thread
$(E_{i-n+1},\cdots, E_i)$. It cannot have loops at the vertex $i$
corresponding to $E_i$ because this would correspond to an
irreducible morphism $E_{i-kn}\to E_{i}$. Similarly a 2-cycle at
$i$ would correspond to irreducible morphisms $E_{i-kn}\to E_j\to
E_{i}$. By periodicity of the helix this would imply that
$-kd+1=\phi(E_j)=-1$ which is impossible for $d\geq 3$.
\end{pf}


The analogue of Proposition \ref{blbl} relating mutation to tilting is

\begin{prop}
\label{w} Let $Z$ be a smooth Fano variety of dimension $d-1$  and
suppose $\H$ is a geometric helix in $\D(Z)$ of type $(n,d)$.
Suppose there is a height function $\phi\colon \H\to\Z$ for an
object $E\in \H$ and write $\sigma_0(\H,\phi)=(\H',\phi')$. Then
the algebra $B'=B(\H')$ is the left tilt of the algebra $B=B(\H)$
at the vertex corresponding to the object $E$.
\end{prop}

\begin{pf}
Take a thread $\E\subset \H$ containing $\E_0=\{E\}$ and
$\E_{-1}$. The mutated helix $\H'$ is generated by the mutated
thread $\E'=\sigma_0(\E)$. Let $Q$ and $Q'$ be the quivers
underlying $B=B(\H)$ and $B'=B(\H')$ respectively. The vertices of
these quivers are in natural bijection with the elements of the
collections $\E$ and $\E'$ respectively. Set
$E'=L_{\E_{-1}}(E)[-1]$. Since $\E\setminus \{E\}=\E'\setminus\{
E'\}$  there is an obvious choice of bijection $\psi\colon Q_0\to
Q'_0$. Consider the diagram of equivalences
\[\begin{CD} \D(B) &@>\Psi>> &\D(B') \\
@V\Phi_{\E}VV && @VV\Phi_{\E'}V \\ \D(Y) &@>\id>> &\D(Y)
\end{CD}\]
Let $i\in Q_0$ be the vertex corresponding to the object $E$.
Applying Lemma \ref{andy} we see that if $j\in Q_0$ and $j\neq i$
then  $\Psi(P_j)=P'_{\psi(j)}$. It remains to  show that
$\Phi_\E(R_i)=\pi^*(E')$ and hence that $\Psi(R_i)=P'_{\psi(i)}$.

Let $A=A(\E)$ be the homomorphism algebra of the thread $\E$.
There is an embedding of  graded algebras $i\colon A\into B$
sending an element of $\Hom_\D(E_i,E_j)$ to the element of $B$
obtained by applying all powers of the functor $S_\D[1-d]$. The
functor $-\tensor_A B$ sends the projective $A$-module $e_i A$ to
the projective $B$-module $e_i B$.
 The existence of this embedding
implies the fact already noted in Section 3 that the quiver $Q_B$
underlying $B$ is obtained from the quiver $Q_A$ underlying $A$ by
adding extra arrows.  The proof of Proposition \ref{plb} implies
that every arrow in $Q_B$ ending at the vertex $i$ actually comes
from an arrow of $Q_A$. It follows that
\[R_i^{B} =R_i^{A} \tensor_A B\]
where $R_i^{A}$ and $R_i^{B}$ are the objects defined in Section
\ref{saville} for the quiver algebras $A$ and $B$ respectively.

Finally,  there is a commuting diagram of functors
\[\begin{CD} \D(A) &@>-\tensor_A B>> &\D(B) \\
@V\Phi_{\E}VV && @VV\Phi_{\E}V \\ \D(Z) &@>\pi^*>> &\D(Y)
\end{CD}\]
Hence the claim follows from Proposition \ref{han} with $k=1$.
\end{pf}

\begin{remark}
\label{aa}
Take assumptions as in Proposition \ref{w}. Write \[(\H^L,\phi^L)=\sigma_0(\H,\phi)\text{ and }(\H^R,\phi^R)=\sigma_0^{-1}(\H,\phi).\]
Then
$B(\H^L)$ is the left tilt of $B(\H)$
at the vertex corresponding to $E$, and $B(\H^R)$ is the right tilt of $B(\H)$
at the vertex corresponding to $E$.

We claim that in the case $d=3$ the helices $\H^L$ and $\H^R$ are the same, up to reindexing, and hence $B(\H^L)$ and $B(\H^R)$ are isomorphic algebras.
Indeed, if we take the thread \[(\E_{-1}, E, \E_1)\subset \H\] then the corresponding threads in $\H^L$ and $\H^R$ are
\[(\L_{\E_{-1}}(E)[-1], \E_{-1},\E_1)\text{ and } (\E_{-1},\E_1, \R_{\E_1}(E)[1]),\]
and the claim follows from Remark \ref{hlx}(b).
\end{remark}


\section{The case of del Pezzo surfaces}

\label{delpezzo} In this section $Z$ is a  del Pezzo surface and
$\D=\D(Z)$. Following work of  Kuleshov and Orlov much
is known about exceptional objects on $Z$. Here we quote their
main results following \cite{gk} and use an argument of Herzog to construct height
functions for strong exceptional collections in $\D$. Combined
with the results of the last section this gives a proof of Theorem
\ref{main} from the introduction. We conclude by discussing some
examples.

\subsection{Exceptional objects on del Pezzo surfaces}

\renewcommand{\lessgtr}{\sim}
Given a torsion-free sheaf $E$ on $Z$ we define the slope
\[\mu(E)=\frac{c_1(E)\cdot (-K_Z)}{r(E)}.\]
If $E$ is a torsion sheaf we set $\mu(E)=+\infty$.  We say that a
torsion-free sheaf $E$ is $\mu$-stable if for all subsheaves $A
\subset E$ one has $\mu(A)<\mu(E)$.

\begin{thm}[Kuleshov, Orlov]
\label{prev} Every exceptional object in $\D$ is a shift of a
sheaf. Moreover every exceptional sheaf $A$ is either a
$\mu$-stable locally-free sheaf or a torsion sheaf of the form
$\O_C(d)$ with $C\subset Z$ an irreducible rational curve
satisfying $C^2=-1$ and $d\in \Z$ is an integer.
\end{thm}

\begin{pf}
This is a combination of Corollary 4.3.2, Theorem 4.3.3 and
Proposition 5.3.3 of \cite{gk}.
\end{pf}

Given an exceptional object $F\in \D$ we write
$f(F)\in \Z$ for the unique integer such that $F\in\Coh(Z)[f(F)]$.
\smallskip

Suppose now that $F_1$ and $F_2$ are two exceptional objects such that $f(F_1)=f(F_2)$. We say $F_1< F_2$ if
\begin{itemize}
\item[(a)] $\mu(F_1)<\mu(F_2)$, or
\item[(b)]   $F_i=\OO_C(d_i)$ and $d_1<d_2$.\end{itemize}
If neither $F_1\leq F_2$ or $F_2\leq F_1$ holds we write
$F_1\lessgtr F_2$. This is the case precisely when $F_1$ and $F_2$
are distinct stable bundles of the same slope, or torsion-sheaves
supported on distinct rational curves.

\begin{thm}
\label{bibl}
For any exceptional pair $(F_1,F_2)$ in $\D$ the complex
$\Hom^\blob_\D(F_1,F_2)$ is concentrated in a single degree. More precisely, if
$f(F_1)=f(F_2)$ then
\begin{align*}
 F_1<F_2 &\implies \Hom^k_\D(F_1,F_2)=0 \text{ unless }k=0, \\
  F_1>F_2 &\implies \Hom^k_\D(F_1,F_2)=0 \text{ unless }k=1,\\
F_1\sim F_2 &\implies \Hom^k_\D(F_1,F_2)=0 \text{ for all
}k.\end{align*}
\end{thm}

\begin{pf}
It is enough to consider the case when $F_1$ and $F_2$ are
sheaves. Consider the triangle \eqref{left}
\[\Hom_{\D}^{\blob}(F_1,F_2)\tensor
F_1\lRa{ev} F_2\lra \L_{F_1} (F_2).\]
By Theorem \ref{prev} the mutated object $L_{F_1}(F_2)$
is of the form $A[n]$ with $A$ a sheaf and $n\in\Z$. Considering the long exact sequence in cohomology it is easy to see that $n\in\{0,1\}$ and $\Hom^\blob_\D(F_1,F_2)$ is concentrated in degrees 0 and 1. If $n=1$ then
$\Hom_{\D}^{\blob}(F_1,F_2)$ must be concentrated in degree 0. If $n=0$ there is a long exact sequence of sheaves
\[0\lra \Hom_\D(F_1,F_2)\tensor F_1\to F_2\lRa{f} \L_{F_1}(F_2)\to  \Hom^1_\D(F_1,F_2)\tensor F_1\to 0.\]
Suppose $\Hom^k(F_1,F_2)$ is nonzero in degrees $k=0$ and $1$ and let $P$ be the cokernel of the map $f$. The sequence
 \[0\to P\lra \L_{F_1}(F_2)\to  \Hom^1_\D(F_1,F_2)\tensor F_1\to 0\]
 gives a nonzero element of $\Ext^1_Z(F_1,P)$. The sequence
 \[0\lra \Hom_\D(F_1,F_2)\tensor F_1\to F_2\lra P\lra 0\]
 then gives a nonzero element of $\Ext^2_Z(F_1,F_1)$ contradicting the fact that $F_1$ is exceptional.
 This proves the first
statement.

Consider now  the second statement. We know from the above that $\Hom^\blob_\D(F_1,F_2)$ is concentrated in degrees 0 and 1. Suppose first that $F_1$ and $F_2$ are stable bundles of slopes $\mu_1$ and $\mu_2$ respectively. If $\mu_1\geq \mu_2$ then
it is standard that $\Hom_\D(F_1,F_2)=0$. On the other hand if $\Hom^1_\D(F_1,F_2)$ is nonzero we have a sequence
\[0\lra F_2\lra \L_{F_1}(F_2)\lra  \Hom^1_\D(F_1,F_2)\tensor F_1\lra 0\]
and since $\L_{F_1}(F_2)$ is a  stable bundle by Theorem \ref{prev}, we can conclude that $\mu_1> \mu_2$. Combining these statements gives the result.

There are three other cases. Suppose first that $F_1$ and $F_2$
are both torsion sheaves, supported on rational curves $C_1$ and
$C_2$. The exceptional pair assumption gives $\chi(F_2,F_1)=0$ and
Riemann-Roch then implies that $C_1\cdot C_2=0$. Since $C_1$ and
$C_2$ are irreducible $(-1)$-curves it follows that they are
disjoint, and  hence $F_1\sim F_2$ and $\Hom^\blob_\D(F_1,F_2)=0$.

Next suppose that $F_1$ is locally-free and $F_2=\OO_C(d)$. Thus
$F_1<F_2$. By Serre duality and the fact that $\K_Z \cdot C<0$ one
has
\[\chi(F_1,F_2)=\chi(F_1^*\tensor \OO_C(d))>\chi(F_1^*\tensor
\omega_Z \tensor \OO_C(d))=\chi(F_2,F_1)=0.\] Hence
$\Hom^k_\D(F_1,F_2)$ must be concentrated in degree 0.
Finally there is the case when $F_1=\OO_C(d)$ and $F_2$ is
locally-free. Then $F_1>F_2$. Arguing as for the previous case one concludes that
$\chi(F_1,F_2)<0$ and hence $\Hom^k_\D(F_1,F_2)$ must be
concentrated in degree 1.
\end{pf}

\subsection{Height functions on del Pezzo surfaces}
Suppose $\E=(E_1,\cdots,E_n)$ is a full strong exceptional
collection in $\D$. In this section we prove that we can reorder
$\E$ so that height functions exist for any object $E\in \E$. The
argument is due to Herzog \cite{he}. As usual  $\F=(F_n,\cdots,
F_1)$ is the dual collection to $\E$.

\begin{lemma}
\label{reorder}
We can reorder the exceptional collection $\E$ (and hence also the dual collection $\F$) so that
the following holds. Suppose $i\geq j$. Then $f(F_i)\geq f(F_j)$, and if equality holds then
either $F_i>F_j$ or $F_i\lessgtr F_j$.
\end{lemma}

\begin{pf}
Consider neighbouring elements $F_{i+1}$ and $F_i$ in $\F$.  If $F_{i+1}$ and $F_i$ are orthogonal  then
so are the corresponding objects $E_i$ and $E_{i+1}$ and we can exchange them in the exceptional collection if necessary.
So let us assume that this is not the case. Write $F_{i+1}=A[m]$ and $F_{i}=B[n]$ with $A$ and $B$
 sheaves. By Theorem \ref{bibl}
\[\Hom_\D(A,B)\neq 0
\text { or } \Hom^1_\D(A,B)\neq 0\] and so $\Hom^k_\D(F_{i+1},F_{i})$ is
non-vanishing in degree $m-n$ or $m-n+1$. Since the objects of $\F$ correspond to simple
modules under the equivalence $\Phi_\E$ one has
\[\Hom^k_\D(F_i,F_j)=0 \text{ for }k<1\]
so it follows that $m-n \geq 0$.

For the second statement, note that  the only other possibility is that
$F_i<F_j$. But by Theorem \ref{bibl} and the assumption that $F_{i+1}$ and $F_i$ are not orthogonal this implies that
$\Hom_\D(F_i,F_j)\neq 0$ which is impossible as before.\end{pf}

Let us reorder our collection $\E$ as in Lemma \ref{reorder} and
fix an object $E\in\E$.
Let $F\in\F$ be the dual object and write
$F=A[p]$ with $A$ a sheaf. Split each subcollection
$f^{-1}(q)\subset\F$ into two subcollections $f^{-1}(q)=R_q\sqcup
L_q$ in such a way that for all $F_j\in f^{-1}(q)$
\[F_j[-q]> A\implies F_j\in R_j  \quad F_j[-q]< A \implies F_j\in L_j.\]
We could take $F\in L_p$ or $F\in R_p$, but for definiteness we choose the first possibility.
We now have a decomposition of $\F$ of the form
\[\F=(\cdots, R_q,L_q,R_{q-1},L_{q-1},\cdots).\]
Of course this induces a decomposition of $\E$ indexed in the opposite direction
\[\E=(\cdots, L_{q-1},R_{q-1},L_{q},R_q,\cdots).\]
There is thus a levelling $\phi\colon\E\to\Z$ defined by
\[\E_q=R_{q-1}\sqcup L_{q}.\]
It satisfies $\phi(E)=p$.
We can now prove

\begin{lemma}
The levelling $\phi\colon\E\to\Z$ is tilting at level $p$.
\end{lemma}

\begin{pf}
Recall that $E=A[p]$ with $A$ a sheaf and that $\phi(E)=n$.
Take an object $E'\in\E$ with $\phi(E')=p'$ and let $F'\in\F$ be the corresponding dual object.
 Without loss of generality we can assume that $E'$ comes after $E$ in $\E$ and hence that $p'\geq p$.

Set $m=f(F')$ so that $F=A'[m]$ for some sheaf $A'$. Now either $F'\in L_m$ in which case $p'=m$ or $F_j\in R_m$ in which case $p'=m+1$. In the first case
 one has $A'\leq A$ or $A'\lessgtr A$.  Either way $\Hom^\blob_\D(A',A)$ is
concentrated in degree zero, so \[\Hom_\D^k(F',F)=0 \text{ unless }
k=p'-p.\] Similarly, if $F_j\in R_m$ then  $B\geq A$ or $B\lessgtr A$ and either way
$\Hom^\blob_\D(A',A)$ is concentrated in degree 1, so again
\[\Hom_\D^k(F',F)=0 \text{ unless }
k=p'-p.\]Thus $E'$ is $p'-p$ related to $E$.  \end{pf}

\subsection{Rolled up helix algebras for del Pezzo surafces}Putting everything together we have now proved

\begin{thm}
Let $\H$ be a geometric helix on a del Pezzo surface. Then the rolled-up helix algebra $B=B(\H)$
is a graded \CY quiver algebra which is noetherian and  finite over its centre. The underlying quiver of $B$ has no loops or 2-cycles. For any vertex $i$ of $Q$ there is another
geometric helix $\H'$ on $Z$ such that the algebra $B(\H')$ is the (left or right) tilt of $B(\H)$ at the vertex $i$.
\end{thm}

\begin{pf}
Suppose $\H$ is a geometric helix on $Z$ and $E\in \H$. By Lemma
\ref{pf} and the existence of height functions for exceptional
collections proved in the last section, we know that we can
reorder $\H$ so that a height function $\phi$ exists for $E\in \H$. Of
course reordering does not affect the underlying rolled-up heix
algebra. The claimed properties of the algebra $B(\H)$ then
follow from Theorem \ref{helequiv}
 and Proposition \ref{plb}. The statement about tilting follows from
 Proposition \ref{w} and Remark \ref{aa}.
\end{pf}

As explained in the introduction, the quivers arising via the tilting process can now be completely understood by the cluster mutation rule. It remains to give some examples of geometric helices on del Pezzo surfaces. Recall that any such surface is either $\PP^1\times \PP^1$ or the blow-up of $\PP^2$ at $0\leq m\leq 8$ points.
We have already given such examples in the case $Z=\PP^2$ and $Z=\PP^1\times\PP^1$.

\begin{example}
\label{thisistheend}
\begin{itemize}
\item[(a)]
On the del Pezzo $Z$ which is $\PP^2$ blown up at one point the exceptional collection
\[\big(\O,\O(h-e),\O(h),\O(2h-e)\big)\]
generates a geometric helix of type $(4,3)$. Here $h$ is the strict transform of a line in $\PP^2$ and $e$ is the exceptional divisor. The canonical bundle is $\O(-3h+e)$.\\

\item[(b)]On the del Pezzo $Z$ which is $\PP^2$ blown up at two points the exceptional collection
\[\big(\O,\O(h-e_1),\O(h-e_2), \O(h),\O(2h-e_1-e_2)\big)\]
generates a geometric helix of type $(5,3)$. Here again $h$ is the strict transform of a line in $\PP^2$,  and $e_1$ and $e_2$ are the exceptional divisors. The canonical bundle is $\O(-3h+e_1+e_2)$. \\

\item[(c)] On a del Pezzo surface $Z$ which is the blow up of $\PP^2$ in $3\leq m\leq 8$ points,
 Karpov and Nogin \cite[Proposition 4.2]{kn} constructed 3-block exceptional collections of sheaves on $Z$. By Proposition \ref{below} these generate geometric helices
 of type $(m+3,3)$.
\end{itemize}
\end{example}


\begin{appendix}
\section{Quiver algebras and tilting}

In this section we give sketch proofs of some simple and
well-known results about quivers for which we could find no
suitable reference. We take notation as in the introduction. In
particular a quiver algebra is one of the form
\[A=A(Q,I)=\CQ/I\]
with $I\subset \CQ_{\geq 2}$, and the augmentation ideal $A_+\subset A$ is spanned by paths of length $\geq 1$.

We use the convention that paths compose on the left, i.e. that $a_1\cdot a_2=0$ unless the target of the arrow  $a_2$ is the source of the arrow $a_1$.
Thus the space of paths (modulo relations) from vertex $i$ to vertex $j$ is $e_j A e_i$.
Since we are considering right-modules this means that a module $M$ determines vector spaces $V_i=Me_i$ for each vertex $i\in Q_0$ and linear maps
$V_j\to V_i$ for each arrow from vertex $i$ to vertex $j$.

%

%
%
\smallskip
The question of which augmented algebras  can be presented as quiver algebras seems to be a tricky one. All we shall need is the following simple result.

\begin{lemma}
\label{simon} Suppose $A=\bigoplus_{n\geq 0} A_n$ is a
finitely-generated graded algebra such that $A_0=\S$ is a
finite-dimensional semisimple algebra. Let $A_+$ denote the augmentation ideal $\bigoplus_{n>0} A_n$.
Then as an augmented algebra $(A,A_+)$ is isomorphic to a a quiver
algebra $A(Q,I)$. Moreover the quiver $Q$ is uniquely determined by the pair  $(A,A_+)$.  \end{lemma}

\begin{pf}
Each graded piece $A_k$ is a finite-dimensional $\S$-bimodule via
left and right multiplication by elements of $A_0=\S$. For each
$k\geq 1$ let $V_k$ be the cokernel of the map of $\S$-bimodules
\[\bigoplus_{0<j<k} A_j\cdot A_{k-j}\lra A_k\]
and choose a splitting $i_k\colon V_k\to A_k$.
Since $A$ is finitely generated we must have $V_k=0$ for $k\gg 0$ and so $V=\bigoplus_k V_k$ is a finite-dimensional $\S$-bimodule and $i=\bigoplus i_k$ is an injective  map of bimodules
$i\colon V\to A$. This induces a map of $\S$-algebras
\[f\colon T_S(V)\to A\]
where $T_S(V)$ is the tensor algebra over $\S$ of the bimodule
$V$. Note that $T_S(V)$ is an augmented algebra with augmentation ideal spanned by tensors of positive degree. By construction $f$ is surjective, and
has kernel contained in the square of the augmentation.

Choose a basis of orthogonal idempotents $(e_1,\cdots,e_n)$ in $\S$ and let $Q$ be the quiver with vertices $\{1,\cdots,n\}$ and $\dim_\C e_i V e_j$ arrows from vertex $i$ to vertex $j$. Then it is easy to see that as an augmented algebra $T_S(V)$ is isomorphic to the path algebra $\CQ$. It follows that $A$ is a quiver algebra. Uniqueness follows from equation
\eqref{nearlydone} from the introduction.
\end{pf}

We shall need the existence of minimal projective resolutions as
in the following Lemma.

\begin{lemma}
\label{min} Suppose $A=A(Q,I)$ is a quiver algebra. For each
vertex $i\in Q_0$ let $S_i$ be the corresponding simple module.
Then there is a projective resolution of the form
\[\cdots \lra \bigoplus_{j\in Q_0} P_j ^{\oplus d_{j,i
}^k} \lra \cdots  \lra \bigoplus_{j\in Q_0} P_j^{\oplus d_{j,i}^1}
\lra P_i \lra S_i\lra 0\] where $d_{j,i}^p=\dim_\C
\Ext^p_{A}(S_i,S_j).$
\end{lemma}

\begin{pf}
An alternative way to state this is that one can construct a projective resolution such that when one applies the functor $\Hom_A(-,S_j)$ all maps
become zero.
One can build such a resolution step-by-step. All one needs to know is
that for any finitely-generated module
$M$ there is a projective module \[P=\bigoplus_{k\in Q_0}
P_k^{n_k}\] and a surjection $f\colon P\to M$ such that for each
$j$ the induced map
\[f^*\colon \Hom_A(M,S_j) \to \Hom_A(P,S_j)\] is a surjection. To
prove this write $M=\bigoplus_{i\in Q_0} M_i$ where $M_i= M e_i$ and take
elements $m_i\in M_i$ such that the images in $M/A_+ M$ form a
basis. Each element $m_i$ defines a map $P_i\to M$, and the
corresponding map $\bigoplus P_i \to M$ has the required property.
\end{pf}

\begin{remark}
\label{rich}If $Q$ has no oriented cycles then it is possible to
order the vertices of the quiver so that $\Hom_A(P_i,P_j)=0$
unless $i<j$. It follows that the projective resolutions above
must be finite, and hence $A$ has finite global dimension.
\end{remark}

Finally we prove the claim made in Remark \ref{bandy}. Suppose
$A=A(Q,I)$ and $A'=A(Q',I')$ are quiver algebras related by a tilt
at the vertex $i$ as in Definition \ref{tommy}. Assume  that their
underlying quivers have no loops.

\begin{lemma}
\label{golden} Define objects $U_j\in D_{\fin}(A)$ by the
relation $\Psi(U_j)=S_{\psi(j)}$.  Then
 $U_i=S_i[-1]$, whereas for $j\neq i$ the object $U_j$ is the universal extension
\begin{equation}  0\lra S_j\lra U_j\lra
\Ext^1_{A}(S_i,S_j)\tensor {S_i} \lra 0.\end{equation}
\end{lemma}

\begin{pf}
Let $U_j$ be the given objects; we will prove that
$\Psi(U_j)=S_{\psi(j)}$. Define an object
\begin{equation} T = R_i\oplus \bigoplus_{i\neq j\in Q_0}
P_j\in\D(A).\end{equation} Then $\Psi(T)=A'$ and to prove the claim we
must check that $\Hom_A(T,U_j)=\C$ for each vertex $j\in Q_0$. The
only tricky thing is to show that $\Hom_A(R_i,U_j)=0$ for $j\neq
i$. Let $V_{ij}$ be the $n_{ij}$-dimensional vector space spanned
by the arrows from $i$ to $j$. Then since we are considering right-modules
$\Ext^1_A(S_i,S_j)=V_{ji}^*$. We must show that the canonical map
$P_j\tensor V_{ji} \to P_i$ induces an isomorphism
\begin{equation}\label{really}V_{ji}^*=\Hom_A(P_i,U_j)\lra \Hom_A(P_j\tensor
V_{ji},U_j)=V_{ji}^* .\end{equation}

Viewed as a representation of the
quiver $Q$, the object $U_j$ can be represented by associating
the vector space $V_{ji}^*$ to the vertex $i$, the one-dimensional
vector space $\C$ to the vertex $j$, and the tautological  linear
map $V_{ji}^* \to \C$ to each arrow $a\in V_{ji}$. On the other
hand $P_i=e_i A$ associates to each vertex $j$ the space of paths in $Q$
(modulo relations) from $j$ to $i$. For each element $\xi\in
V_{ji}^*$ there is a map $P_i\to U_j$ defined by sending the lazy
path at $i$ to $\xi\in V_{ji}^*$, sending an arrow $a\in V_{ji}$
to $\xi(a)\in \C$, and sending all other paths to zero. Composing
with the canonical map $P_j\tensor V_{ji} \to P_i$ gives a nonzero
map for all nonzero $\xi$. This shows that \eqref{really} is
injective, and hence an isomorphism.
\end{pf}

\section{Semi-orthogonal decomposition}
\label{semi} Here we recall Bondal's categorical approach to
mutation functors. For more details we refer the reader to \cite{bo}. Throughout $\D$
denotes an arbitrary $\C$-linear triangulated category of finite type.

Suppose $\A\subset \D$ is a full subcategory.  The \emph{right
orthogonal subcategory} to $\A$ is
\[ \A^{\perp} =\{X\in \D : \Hom_{\D}^{\blob} (A,X)=0 \text{ for }
A\in\A\}\subset\D.\] Similarly, the \emph{left orthogonal
subcategory} to $\A$ is
\[{^\perp} \A = \{X\in \D: \Hom_{\D}^{\blob}(X,A)=0 \text{ for } A\in \A\}\subset\D.\]
Both are full triangulated subcategories of $\D$.

\smallskip
A \emph{semi-orthogonal decomposition} of
$\D$ is a pair of full triangulated subcategories $(\A,\B)\subset\D$ such that
\begin{itemize}
\item[(a)] for $A\in\A$ and $B\in\B$ one has $\Hom_{\D}(B,A)=0$,
\item[(b)] for every object $X\in\D$ there is triangle
\[B \lra X\lra A\]
such that $A\in\A$ and $B\in \B$.
\end{itemize}

\smallskip

A full triangulated subcategory $\A\subset\D$ is  \emph{left or
right admissible} if the inclusion functor $\A\into \D$ has a left
or right adjoint respectively.


\begin{prop}[Bondal]
\label{ad}
Suppose $\A,\B\subset\D$ are full, triangulated subcategories closed under isomorphism. Then the following are equivalent
\begin{itemize}
\item[(a)] $(\A,\B)\subset\D$ is a semi-orthogonal decomposition, \smallskip\item[(b)] $\A$ is left
admissible and $\B=^{\perp}\A$,\smallskip \item[(c)] $\B$ is right
admissible and $\A=\B^\perp $.
\end{itemize}
\end{prop}

\begin{pf}
This can be found in \cite{bo}. Here we just sketch the argument.
First assume that $(\A,\B)$ is a semi-orthogonal decomposition. Condition (a) in the definition
of semi-orthogonality implies that for any $X\in\D$ the triangle
 appearing in part (b)  is
unique up to a unique isomorphism. This implies that there are
functors \[p\colon \D\to\A,\quad q\colon\D\to\B,\]
sending an object $X\in\D$ to the objects $A\in\A$ and $B\in\B$
respectively.
It is then easy to see that $p$ is the left adjoint to the inclusion $\A\into\D$ and that
similarly $q$ is the right adjoint to the inclusion $\B\into\D$.

Finally, note that if $X\in{^\perp}\A$ then the map $X\to A$ must be zero. Thus $\id_A:A\to A$ factors via $B[1]\in\B$ and hence is zero by the orthogonality condition. This implies that ${^\perp}\A\subset \B$. The opposite inclusion is immediate from part (a) of the definition.  Similarly $\B^\perp=\A$, so (a) implies (b) and (c).

For the converse we must prove that (a) is implied by either (b) or (c).
Without loss of generality assume (b) so that the inclusion $i\colon\A\to \D$ has a left adjoint $p\colon\D\to\A$.  Taking the cone on
the unit of the adjunction gives for any object $X\in \D$ a
triangle
\[X \lra (i\circ p)(X)  \lra B.\]
Since $i$ is fully faithful $p\circ
i\isom \id_{\A}$ and so
applying $p$ to the above triangle shows that $p(B)=0$. Now
\[\Hom_\D(B,i(A))=\Hom_\A(p(B),A)=0\]
so  $B\in {^\perp}\A=\B$. This proves that $(\A,\B)$ is a semi-orthogonal decomposition.
\end{pf}

A full triangulated subcategory $\A \subset \D$ is called
\emph{admissible} if it both left and right admissible.
Then by Proposition \ref{ad} one has semi-orthogonal decompositions $(\A,{^\perp}\A)$ and $(\A^\perp,\A)$.
Let $p\colon \D\to {^\perp \A}$ be the left adjoint to the inclusion functor $i\colon {^\perp} \A \to \D$,
and let $q\colon
\D\to \A^\perp$ be the right adjoint to the inclusion $j\colon \A^\perp\to \D$.
\[\xymatrix{ {^\perp \A} \ar@/^1pc/[rr]^{i} && \D \ar@/^1pc/[ll]^{p}
\ar@/^1pc/[rr]^{q} && \A^\perp \ar@/^1pc/[ll]^{j}}\] The
composite functors $\L_\A=q\circ i$ and $\R_\A= p\circ j$ are
called the \emph{mutation functors} for the subcategory $\A\subset\D$.

\begin{lemma}
\label{earl}Suppose $X\in {^\perp}\A$ and $Y\in \A^\perp$. Then $Y=\L_\A(X)$ iff there is a triangle
\[A\lra X\lra Y\]
with $A\in\A$. Similarly $X=\R_\A(Y)$ iff there is a triangle
\[X\lra Y\lra A'\]
with $A'\in\A$.
\end{lemma}

\begin{pf}
This follows immediately from the proof of Proposition \ref{ad}.
\end{pf}

Rotating the triangle it is then obvious that the mutation functors $\L_\A$ and $\R_\A$ are mutually-inverse equivalences of categories.

\smallskip

A saturated triangulated category is always admissible in any
enveloping category. In particular, if $\E\subset\D$ is an exceptional collection
then $\langle\E\rangle\subset\D$ is always admissible. It follows immediately from
 Lemma \ref{earl} that if $\E\subset \D$ is an exceptional
collection, then the mutation functors $L_\E$ and $R_\E$ defined
in Section 2 coincide with the mutation functors
$L_{\langle \E\rangle}$ and $L_{\langle \E\rangle}$ defined above.

\section{Exceptional collections}
Here we give the proofs of various simple results on exceptional
collections from Section 2. Assumptions are as in Section 2.1.

\begin{lemma}
\label{full}
Let $\E\subset\D$ be an exceptional collection. Then the following are equivalent
\begin{itemize}
\item[(a)] $\langle\E\rangle=\D$,
\item[(b)] $\E^{\perp}=0$,
\item[(c)] ${^\perp}\E=0$.
\end{itemize}
\end{lemma}

\begin{pf}
Take $X\in\E^\perp$. Then  $\Hom_\D^\blob(E,X)=0$ for all $E\in\E$ and so
$\Hom_\D^\blob(Y,X)=0$ for all $Y\in\langle \E\rangle$. If (a) holds we can take $Y=X$ and so $X=0$. Thus (a) implies (b) and similarly (c).
For the converse note that $\langle\E\rangle$ is saturated and hence is admissible in $\D$. Thus there is a semi-orthogonal decomposition $(\E^{\perp},\langle\E\rangle)$.
If (b) holds then it follows that $\langle\E\rangle=\D$.
\end{pf}

\begin{lemma}
\label{weedoris}
Let $\E=(E_1,\cdots,E_n)$ be a full exceptional collection in $\D$. For any $1\leq k\leq n$ define exceptional collections $\E_{\leq k}=(E_1,\cdots,E_k)$ and $\E_{>k}=(E_{k+1},\cdots,E_n)$. Then for each $k$ there is a semi-orthogonal decomposition
 \[\langle\E\rangle=(\langle \E_{\leq k}\rangle,\langle\E_{>k}\rangle).\]
\end{lemma}

\begin{pf}
The subcategory $\langle\E_{>k}\rangle$ is saturated and hence admissible. There is thus a semi-orthogonal decomposition
$(\E_{>k}^\perp,\langle\E_{>k}\rangle)$. Clearly $\E_{\leq k}$ is an exceptional collection in $\E_{>k}^\perp$ and since $\E$ is full in $\D$ it follows
from Lemma \ref{full} that $\E_{\leq k}$ is full in $\E_{>k}^\perp$. This gives the result.
\end{pf}

\begin{lemma}
\label{grot}
If $\D$ contains a full exceptional collection of length $n$ then $K(\D)\isom \Z^{\oplus n}$.
\end{lemma}

\begin{pf}
If $E\in\D$ is an exceptional object then there is an equivalence $D(\C)\to \langle E\rangle$ sending $\C$ to $E$.
In particular $K(\langle E\rangle)=\Z$ is generated by the class $[E]$.
This can be seen as a special case of Theorem \ref{RB} but is also easy to check directly. Now suppose $\E=(E_1,\cdots,E_n)$ is a full exceptional collection in $\D$.
Applying Lemma \ref{weedoris} repeatedly it follows that $K(\D)$ is spanned by the classes $[E_i]$.
Recall that the Euler form
\[\chi(X,Y)=\big.\sum_{i\in\Z} (-1)^i \dim_\C \Hom_\D^i(X,Y)\]
descends to $K(\D)$. Since $\chi(E_i,E_j)=0$ for $i>j$ and $\chi(E_i,E_i)=1$ it follows that the classes $[E_i]\in K(\D)$
are linearly independent in $K(\D)$ and hence form a basis.
\end{pf}

\medskip

\noindent {Proof of Lemma \ref{serre}.}
Let $q,r\colon \D\to \E^\perp$ be the left and right adjoints to
the inclusion functor $j\colon \E^\perp\to \D$ respectively.  The
first claim is that \[S_{\E^\perp}\circ q\isom r\circ S_{\D}.\]
This follows from the definition of the Serre functor. Indeed,
writing $\CC=\E^{\perp}$ one has
\begin{align*}\Hom_{\CC}(X, S_\CC (
q(Y)))&=\Hom_{\CC}(q(Y),X)^*=\Hom_\D(Y,j(X))^*
\\
&=\Hom_\D(j(X),S_\D(Y))=\Hom_{\CC}(X,r(S_\D(Y)))\end{align*}
 for
$X\in \CC$ and $Y\in \D$. Let $i\colon {^\perp \E} \to \D$ be the
inclusion functor. Using the description $\L_\E=q\circ i$ from the last section, the equivalence
\eqref{ffs} can be rewritten as
\[S_\D|_{{^\perp}\E}\isom p\circ S_\D \circ i\isom  S_{\E^\perp}\circ q \circ i\isom S_{\E^\perp}\circ L_\E.\]
This gives the result. \qed

\medskip

\noindent {Proof of Lemma \ref{dual}.} Write $\E_{\leq j}$ for the
subcollection $(E_1,\cdots,E_j)\subset \E$ and put
$\E_{<j}=\E_{\leq j-1}$. Then $F_j=\L_{\E_{<j}}(E_j)$ and is
therefore an object of the subcategory
\[\E_{<j}^\perp\subset \langle\E_{\leq j}\rangle.\]
It follows immediately that $\F$ is an exceptional collection, and
that \[\Hom_\D^k(E_i,F_j)=0\text{ for }i\neq j.\] It is also easy to see
that the collection $\F$ is full. The triangle of Lemma \ref{earl}
takes the form
\[ Y\lra  E_j \lra F_j \] for some $Y\in \E_{<j}$. Since $E_i\in
{^\perp}\E_{<i}$, applying the functor $\Hom_\D(E_i,-)$ shows that
\[\Hom_{\D}^\blob(E_i,F_i)=\C.\] This completes the proof of the
first part of the Lemma.

For the uniqueness statement note first that by assumption $F_j\in
\E_{<j}^{\perp}$.  Let $Y$ be the cone on a nonzero map $E_j\to
F_j$ fitting into a triangle \[  Y\lra E_j \lra F_j \] as before.
Applying the functor $\Hom_\D(E_j,-)$ shows that $\Hom(E_j,Y)=0$
and therefore by Lemma \ref{weedoris} \[Y\in \E_{\geq
j}^{\perp}=\langle\E_{< j}\rangle.\] Applying Lemma \ref{earl} it
follows that $F_j=\L_{\E_{<j}}(E_j)$.

\end{appendix}

\end{document}